\documentclass[12pt]{article}
\usepackage{amssymb,amsmath}

\input epsf.tex

\title{Homological mirror symmetry and 
torus fibrations}
\author {Maxim Kontsevich and Yan Soibelman}

\begin{document}
\maketitle

\newtheorem{thm}{Theorem}
\newtheorem{lmm}{Lemma}
\newtheorem{dfn}{Definition}
\newtheorem{rmk}{Remark}
\newtheorem{prp}{Proposition}
\newtheorem{conj}{Conjecture}
\newtheorem{exa}{Example}
\newtheorem{cor}{Corollary}
\newtheorem{que}{Question}
\newtheorem{ack}{Acknowledgements}
\newcommand{\K}{{\bf k}}
\newcommand{\C}{{\bf C}}
\newcommand{\R}{{\bf R}}
\newcommand{\N}{{\bf N}}
\newcommand{\Z}{{\bf Z}}
\newcommand{\Q}{{\bf Q}}
\newcommand{\G}{\Gamma}
\newcommand{\A}{A_{\infty}}

\newcommand{\ihom}{\underline{\Hom}}
\newcommand{\ra}{\longrightarrow}
\newcommand{\epi}{\twoheadrightarrow}
\newcommand{\mono}{\hookrightarrow}

\newcommand{\epp}{\varepsilon}

\section{Introduction}

\subsection{Homological mirror symmetry and degenerations}

Mathematically mirror symmetry can be interpreted in many ways.
In this paper we will make a bridge between two approaches: the 
homological mirror symmetry ([Ko]) and the duality between torus
fibrations (a version of Strominger-Yau-Zaslow conjecture, see [SYZ]).

The mirror symmetry is a duality between Calabi-Yau manifolds, 
i.e. complex manifolds
which carry a K\"ahler metric with vanishing Ricci curvarure.
In fact, it is rather duality not  between individual
 manifolds, but between manifolds in certain ``degenerating'' families
 (``large complex structure limit'' and ``large symplectic structure limit'').
 In this paper we propose a
 differential-geometric model 
 of this degeneration. In particular, we conjecture that
 in the limit both dual manifolds $X$ and $X^{\vee}$ become  fiber bundles
 with toroidal fibers over the {\it same} base $\overline{Y}$
 (see Section 3). Metric space $\overline{Y}$ is a compactification with
some mild singularities of a (real) Riemannian manifold  $Y$ whose dimension is
half of the dimension of $X$ and $X^{\vee}$.    Also, the manifold $Y$ carries
a rich geometric structure,    including certain ``combinatorial'' data
(so-called integral    affine structure).
  This picture is partially
 motivated by the classical theory of collapsing 
 Riemannian manifolds developed by M. Gromov and others
 (see for ex. [CG]). Another origin of our geometric conjectures is the 
 [SYZ] version of mirror duality. We recast it in somewhat different
 terms in Section 2, devoted to the moduli space of conformal
 field theories and its natural compactification. In a recent preprint
  [GW] similar differential-geometric conjectures were suggested and
   verified in the case of degenerating K3-surfaces.

 The Homological Mirror Conjecture proposed in [Ko] is
 a statement about equivalence
of two $\A$-categories:
 the (derived) category of coherent sheaves on
a Calabi-Yau manifold $X$
and the Fukaya category of the dual Calabi-Yau manifold $X^{\vee}$.
The former is defined in holomorphic (or algebraic) terms, the latter is
defined in terms of symplectic geometry.

  We apply the ideas of the theory of collapsing Riemannian manifolds
  to the Homological Mirror 
  Conjecture. Let us
   call {\it differential-geometric} this model for degenerating
  Calabi-Yau manifolds. It gives a clear picture for
  the degeneration of the Fukaya category.
  The Fukaya category $F(X^{\vee}, \omega^{\vee})$
 of a symplectic manifold
 $(X^{\vee}, \omega^{\vee})$, with $[\omega^{\vee}]\in H^2(X^{\vee},{\Z})$,
   and its degeneration are defined
  as $\A$-categories over the field of Laurent formal power series
 ${\C}((q))$. The parameter $q$ enters in the story
 when one writes higher compositions (Massey products),
 which have expressions $q^{\int_{\beta} \omega^{\vee}}$ as coefficients, 
 $\beta \in H_2(X^{\vee},{\Z})$.
 We can set $q=exp(-1/\varepsilon), \varepsilon \to 0$, where the parameter
 $\varepsilon$ corresponds to the rescaling of the symplectic form:
 $\omega^{\vee}\mapsto \omega^{\vee}/\varepsilon$.
 If $[\omega^{\vee}]$ does not belong to $H^2(X,{\Z})$, one can work
over the field 
${\C}_{\varepsilon}:=\{\sum_{i\ge 0}a_ie^{-\lambda_i/\varepsilon}|
\,a_i\in {\C}, \lambda_i\in {\R}, \lambda_i\to +\infty\}$.

 In the case of torus fibrations,
  a full subcategory of the limiting Fukaya category
 can be described in terms of the Morse theory on the base
 of the torus fibration. The higher products giving the
 $\A$-structure can be written as sums over sets of 
 planar trees. In the case of cotangent bundles 
 (instead of torus fibrations)
 this description was proposed earlier by Fukaya and Oh (see [FO]).

 On the holomorphic
 side of mirror symmetry, the degeneration of the dual family
 $X_q$ is described
 in non-archimedean terms: we have a Calabi-Yau manifold ${\cal X}_{mer}$
 over the field ${\C}_q^{mer}$ of germs at $q=0$ of meromorphic 
 functions.
 Changing scalars, we get a Calabi-Yau manifold ${\cal X}_{form}$
 over the local field of Laurent series ${\C}((q))$. Let us call {\it analytic}
 this degeneration picture. 
 There is a description of a class of
 algebraic Calabi-Yau manifolds (over arbitrary
  local  fields, complete with respect to discrete valuations)
  in terms of real $C^{\infty}$-manifolds with integral affine structures.
  We 
 expect that differential-geometric and (non-archimedean) analytic
  pictures of the degeneration
 are equivalent. This equivalence reflects two different ways of
 looking at Calabi-Yau manifolds: differential-geometric (via K\"ahler metrics
 with vanishing Ricci curvature) and algebro-geometric (via 
 smooth projective varieties with vanishing canonical class).

The Homological Mirror Conjecture
says that the Fukaya category
$F(X^{\vee}, \omega^{\vee})$ is equivalent (as an
$\A$-category over ${\C}((q))\,$) to the derived category
of coherent sheaves $D^b({\cal X}_{form})$. We expect that it 
implies well-known numerical predictions for the number of rational curves
on a Calabi-Yau manifold (genus zero Gromov-Witten invariants).

Using our conjectures about the collapse of Calabi-Yau manifolds, we offer
a general approach to the
proof of Homological Mirror Conjecture. We apply it
in the case when the torus fibration has no singularities.
This happens in the case of abelian varieties. 
In general, one should investigate the input
of singularities of the base of torus fibration. 
In the same vein, we discuss the relationship between
 Riemannian manifolds with intergral affine structures and 
 varieties over non-archimedean fields
   in the simplest case of flat tori and abelian varieties.
 The general case will be discussed elsewhere.

It should be clear from the above discussion that the non-archimedean
analysis plays an important role in the formulation and the proof of
Homological Mirror Conjecture.
Analytic picture of the degeneration
seems to be related to the theory of rigid analytic spaces in the version
of Berkovich (see [Be]). In particular, there is a striking similarity
between the base of torus fibration and a certain
canonically defined subset (see 3.3) of the Berkovich spectrum
 of an algebraic Calabi-Yau manifold over a local field. 
This subject definitely deserves further investigation.

\subsection{Content of the paper}

In Section 2 we discuss motivations from the Conformal Field Theory.
 In Section 3 we 
formulate the conjectures about analytic and
differential-geometric pictures of 
the large complex structure limit. In Section 4 we describe a general
framework of $\A$-pre-categories adapted to the transversality
problem in the definition of the Fukaya category. 
Section 5 is devoted to the Fukaya category
and its degeneration. The reader will notice an advantage of working over the 
field of Laurent power series: one can consider {\it all} local systems
over Lagrangian submanifolds (in the conventional approach unitarity
of the holonomy is required). Section 6 is devoted to the $\A$-category 
of smooth
functions introduced by Fukaya (and then studied by Fukaya and Oh in [FuO]).
We prove that this $\A$-category has a very simple de Rham model.
 This part of the paper can be read independently
of the rest. On the other hand, the technique and the general scheme
of the proof will be used later
in Section 8, devoted to the Homological
Mirror Conjecture. One important technical tool is an explicit $\A$-structure
on a subcomplex of a differential-graded algebra (see [GS], [Me]).
We restate the formulas from [Me] in term of sums
over a set of planar trees. Our proof of the equivalence of
Morse and de Rham $\A$-categories uses the approach to Morse theory from
[HL]. Section 7 is devoted to the analytic side
of the Homological Mirror Conjecture. We assign a rigid analytic space
to the class of torus fibrations discussed in Section 3.
We also construct a mirror symmetry functor for torus fibrations
in terms of non-archimedean geometry. The use of
non-archimedean analysis allows us to avoid problems with
convergence of series in the definition of the Fukaya category.
In Section 8 we construct
the $\A$-pre-category which is equivalent to a 
full $\A$-subcategory of the $\A$-version of the
derived category of
coherent sheaves on a Calabi-Yau manifold over ${\C}_{\varepsilon}$.
We prove that this category is  
equivalent to an $\A$-subcategory 
of the Fukaya category of the mirror dual torus fibration. 
In Appendix (Section 9) we describe the analogs of our constructions
in the case of complex geometry.

{\it Acknowledgements}. Research of Y.S. was partially supported
by IHES, the fellowship from the Clay Mathematics Institute and
CRDF grant UM-2091. He also thanks the IHES for hospitality and
excellent research conditions.

\section{ Degenerations of unitary Conformal Field Theories} 

In this section we will explain physical motivations for
our picture of mirror symmetry. We assume that the reader is
familiar to some extent with the basic notions of Conformal
Field Theory. For example, the lectures [Gaw] contain
most of what we need.

Unitary Conformal Field
Theory (abbreviated by CFT below ) is well-defined mathematically.
It is described by the following data:

1) A real number $c\ge 0$ called central charge.

2) A bi-graded pre-Hilbert {\it space of states} 
$H=\oplus_{p,q\in {\bf R}_{\ge 0}}H^{p,q},p-q\in {\bf Z}$ such that 
$dim(\oplus_{p+q\le E} H^{p,q})$ is 
finite for every $E\in {\bf R}_{\ge 0}$. Equivalently, there is an
action of the Lie group ${\C}^{\ast}$ on $H$, so that
$z\in{\C}^{\ast} $ acts on $H^{p,q}$ 
as $z^p\bar{z}^q:=(z\bar{z})^p\bar{z}^{q-p}$.

3) An action of the product of Virasoro  and anti-Virasoro Lie algebras
$Vir\times \overline{Vir}$ (with the same central charge $c$)
on $H$, so that the space $H^{p,q}$ 
is an eigenspace of
the generator $L_0$ (resp. $\overline{L}_0$) with the eigenvalue $p$
(resp. $q$).
 
4) The space $H$ carries some additional 
structures derived from the operator product expansion (OPE).
The OPE is described by a 
linear map $H\otimes H\to H\widehat{\otimes}{\C}\{ {z,\bar{z}}\}$.
Here ${\C}\{ {z,\bar{z}}\}$ is the topological ring of formal
power series $f=\sum_{p,q}c_{p,q}z^{p}\bar{z}^{q}$ where $ c_{p,q}\in {\C},
\,p,q\to +\infty,\,\, p,q \in {\R},\, p-q\in {\Z}$.
The OPE satisfies a list axioms, which we are not going to recall here
(see [Gaw]).

Let $\phi\in H^{p,q}$. Then the number $p+q$ is called the 
{\it conformal dimension} of $\phi$
(or the {\it energy}), and $p-q$ is called the {\it spin}
of $\phi$. Notice that, since the spin of $\phi$
is an integer number, the condition $p+q<1$ implies $p=q$.

The central charge $c$ can be described by the formula
$dim(\oplus_{p+q\le E}H^{p,q})=exp(\sqrt{4/3\pi^2cE(1+o(1))}$ as
$E\to +\infty$. It is expected that all possible central charges
form a countable well-ordered subset of ${\bf Q}_{\ge 0}\subset {\R}_{\ge 0}$.
If $H^{0,0}$ is a one-dimensional vector space, the corresponding CFT is called
irreducible. A general CFT
is a sum of irreducible ones. The {\it trivial} CFT has $H=H^{0,0}=\C$
and it is the unique 
irreducible unitary CFT 
with $c=0$.

\begin{rmk} Geometric considerations of this paper are related 
to $N=2$ Superconformal Field Theories (SCFT). There is a version
of the above data and axioms for SCFT. In particular, each
$H^{p,q}$ is a hermitian super vector space. There is an action
of the super extension of the product
of Virasoro and anti-Virasoro algebra on $H$. 
In the discussion of the moduli
spaces below we will not distinguish between CFTs and SCFTs, because
except of some minor details, main conclusions are true in both cases.
\end{rmk}

\subsection{Moduli space of Conformal Field Theories}

For a given CFT one can consider its group of symmetries
(i.e. automorphisms of the space $H=\oplus_{p,q} H^{p,q}$ 
preserving all the structures).
It is expected that the group of symmetries is a compact Lie
group of dimension less or equal than $dim\,H^{1,0}$.

Let us fix $c_0\ge 0$ and $E_{min}>0$, and consider the moduli space
${\cal M}_{c\le c_0}^{E_{min}}$ of
all irreducible CFTs with the central charge $c\le c_0$ and
$$min\{p+q>0| H^{p,q}\ne 0\}\ge E_{min}$$ 
It is expected that ${\cal M}_{c\le c_0}^{E_{min}}$
is a compact real analytic stack  of finite local dimension.
The dimension of the base of the minimal
versal deformation of a given CFT is less or equal than $dim\,H^{1,1}$.
We define
${\cal M}_{c\le c_0}=\cup_{E_{min}>0}{\cal M}_{c\le c_0}^{E_{min}}$.
We would like to compactify this stack by adding boundary components
corresponding to certain asymptotic
descriptions of the theories with $E_{min}\to 0$.
The compactified space
is expected to be a compact stack $\overline{{\cal M}}_{c\le c_0}$.
In what follows we will loosely use the word ``space'' instead of
the word ``stack''.

\begin{rmk}
There are basically only two classes of rigorously defined CFTs: the rational 
theories (RCFT) and the
lattice CFTs. Considerations of this paper correspond
to the case of sigma models which produce neither of these. 
The description of sigma models as path integrals corresponding
 to certain Lagrangians did not give yet
a mathematically satisfactory construction.
As we will explain below, there is an alternative way to speak about 
sigma models
in terms of degenerations of CFTs.

\end{rmk}

\subsection{Physical picture of a simple collapse}

In order to compactify ${\cal M}_{c\le c_0}$
we consider degenerations of CFTs as $E_{min}\to 0$.
A degeneration is given by a one-parameter (discrete or continuous)
 family $H_{\varepsilon}, \varepsilon\to 0$ 
of bi-graded spaces as above,
where $(p,q)=(p(\varepsilon),q(\varepsilon))$. These
spaces are equipped with OPEs.
 The subspace of fields with
conformal dimensions  vanishing as $\varepsilon\to 0$
gives rise to a commutative algebra $H^{small}=\oplus_{p(\varepsilon)\ll 1}
H_{\varepsilon}^{p(\varepsilon),p(\varepsilon)}$
(the algebra structure is given by the leading terms
in OPEs).  
The spectrum $X$ of $H^{small}$ is expected to be a compact space 
(``manifold with
singularities'')  such that $dim\,X\le c_0$.
It follows from the conformal invariance and the OPE, that 
the grading of $H^{small}$ (rescaled as $\varepsilon\to 0$)
is given by the eigenvalues of a
second order differential operator defined on the smooth part of $X$.
The operator has positive eigenvalues and is determined
up to multiplication by a scalar. This implies that the smooth part of $X$ 
carries a metric $g_X$, which is
also defined up to multiplication by a scalar. Other terms in OPEs give rise
to additional differential-geometric structures on $X$.
 
 Thus, as a first approximation to the real picture, we assume
 the following description of a ``simple collapse'' of a family of CFTs.
The degeneration of the family is described by the point of the boundary of
$\overline{{\cal M}}_{c\le c_0}$ which is a
triple $(X, {\R}_+^{\ast}\cdot g_X,\phi_X)$, where the 
 metric $g_X$  is defined up to
a positive scalar factor, and $\phi_X:X\to {\cal M}_{c\le c_0-dim\,X}$
is a map. One can have some extra conditions on the data.
For example, the metric $g_X$ can satisfy the Einstein equation.

Although the scalar factor for the metric is arbitrary, one should
imagine that the curvature of $g_X$ is ``small'', and the injectivity 
radius of $g_X$ is ``large''.
The map $\phi_X$ appears naturally from the point
of view of the simple collapse of CFTs described above.
Indeed, in the limit $\varepsilon\to 0$, the space $H_{\varepsilon}$
becomes an $H^{small}$-module. It can be thought of as a space
of sections of an infinite-dimensional vector bundle
$W\to X$. One can argue that fibers of $W$ generically are
spaces of states of CFTs with central charges less or equal
than $c_0-dim\,X$.
 This is encoded in the map $\phi_X$. In the case when CFTs
from $\phi_X(X)$ have non-trivial symmetry groups, one expects
a kind of a gauge theory on $X$ as well.

Purely bosonic sigma-models correspond the case when 
$c_0=c(\varepsilon)=dim\,X$
and the residual theories (CFTs in the image of $\phi_X$) are all trivial.
The target space $X$ in this case should carry a Ricci flat metric.
In the supersymmetric case the target space $X$ is a Calabi-Yau manifold,
and the residual bundle of CFTs is a bundle of free fermion theories.

\begin{rmk}
We expect that all compact Ricci flat manifolds
(with the metric defined up to a constant scalar factor) appear as target
spaces of degenerating CFTs. Thus, the construction of the compactification
 of the moduli space of CFTs should include as a part a compactification
 of the moduli spaces of Einstein manifolds.  Notice that
 in differential
  geometry there is a fundamental result of Gromov (see [G]) about the 
  precompactness of  
  the moduli 
  space of pointed connected complete Riemannian manifolds of a given dimension,
   with the 
  Ricci curvature bounded from below. 
   There is a deep relationship between the
    compactification of the moduli space of CFTs
     and the Gromov's compactification.
   It seems that one can deduce from certian physical arguments
 that
 all target spaces appearing as
limits of CFTs     have non-negative Ricci curvature.      
\end{rmk}    

\subsection{Multiple collapse and the structure of the boundary}

In terms of the Virasoro operator $L_0$ the
 collapse is described by a subset (cluster) $S_1$ 
in the set of eigenvalues of $L_0$
which approach to zero ``with the same speed'',  as $E_{min}\to 0$.
The next level of the collapse is described by 
another subset $S_2$ of eigenvalues
of $L_0$. Elements of $S_2$ approach to zero ``modulo the first collapse"
(i.e. at the same speed, but ``much slower'' than elements of $S_1$).
One can continue to build a tower of degenerations. 
It leads to an hierarchy of
boundary strata. Namely, if there are further degenerations
of CFTs parametrized by $X$, one gets a fiber bundle
over the space of triples 
$(X,{\R}_+^{\ast}\cdot g_X,\phi_X)$ with the
fiber which is the space of triples of similar sort.
Finally, we obtain the following 
qualitative geometric picture of the boundary 
$\partial \overline {{\cal M}}_{c\le c_0}$. 

A boundary point is given by the following data:

1) A finite tower of maps of compact topological spaces 
$p_i:\overline{X}_i\to
\overline{X}_{i-1}, 0\le i\le k$, $\overline{X}_0=\{pt\}$.

2) A sequence of smooth manifolds 
$({X}_i,g_{{X}_i}), 
0\le i\le k$, such that ${X}_i$ is a dense 
subspace of $\overline{X}_i$, and $dim\,X_i>dim\,X_{i-1}$,
and $p_i$ defines a fiber bundle $p_i:X_i\to X_{i-1}$.

3) Riemannian metrics on the fibers of the restrictions
of $p_i$ to $X_i$, such that the diameter of each fiber is finite.
In particular the diameter of $X_1$ is finite, because it is the
only fiber of the map $p_1:X_1\to \{pt\}$.

4) A map $X_k\to {\cal M}_{c\le c_0-dim\,X_k}$.

The data above are considered up to the natural action
of the group $({\R}_+^{\ast})^{k}$ (it rescales the metrics on fibers).

There are some additional data, like non-linear connections
on the bundles $p_i:X_i\to X_{i-1}$. The set of data should
satisfy some conditions, like differential equations on the metrics.
We cannot formulate this portion of data more precisely in general case.
It will be done below in the case of $N=2$ SCFTs corresponding to
sigma models with Calabi-Yau target spaces.

\subsection{Example: Toroidal models }

Non-supersymmetric toroidal model is described by the
so-called Narain lattice, endowed with some additional data.
More precisely, let us fix the central charge $c=n$ 
which is a positive integer number. What physicists call the Narain lattice $\Gamma^{n,n}$
is a unique unimodular lattice of
rank $2n$ and the signature $(n,n)$. It can be described 
as ${\Z}^{2n}$ equipped with the quadratic form
$Q(x_1,...,x_n,y_1,...,y_n)=\sum_ix_iy_i$.
The moduli space of toroidal CFTs is 

$${\cal M}_{c=n}^{tor}=O(n,n,{\Z})\backslash O(n,n,{\R})/O(n,{\R})
\times O(n,{\R}).$$

Equivalently, it is a quotient  of the open part of
the Grassmannian 
$\{V_+\subset {\R}^{n,n}|\, dim\,V_+=n,Q_{|V}>0\}$
by the action of $O(n,n,{\Z})=Aut(\Gamma^{n,n},Q)$. 
Let $V_-$ be the orthogonal complement to $V_+$.
Then every vector of $\Gamma^{n,n}$
can be uniquely written as $\gamma=\gamma_{+}+\gamma_{-}$,   where
 $\gamma_{\pm}\in V_{\pm}$.
For the corresponding CFT one has
$$\sum_{p,q}dim(\,H^{p,q})z^p\bar{z}^q=
\bigl|\prod_{k\ge 1}(1-z^k)\bigr|^{-2n}\sum_{\gamma\in\Gamma^{n,n}}
z^{Q(\gamma_{+})}\bar{z}^{\,-Q(\gamma_{-})}$$

Let us try to compactify the moduli space ${\cal M}_{c=n}^{tor}$.  Suppose that we have a   one-parameter
 family of toroidal theories such that  $E_{min}(\varepsilon)$ approaches zero. Then for corresponding
vectors in $H_{\varepsilon}$  one gets
$p(\varepsilon)=q(\varepsilon)\to 0$. It implies that
$Q(\gamma(\varepsilon))=0, Q(\gamma_{+}(\varepsilon))\ll 1$. 
It is easy to see that
one can add vectors $\gamma(\varepsilon)$
satisfying these conditions. Thus one gets
a (part of) lattice of the rank less or equal than $n$.
In the case of ``maximal'' simple collapse the rank will be equal to $n$.
One can see that the corresponding points of the boundary give
rise to the following data: $(X,{\R}_+^{\ast}\cdot g_X, \phi_X^{triv};B)$, where
$(X,g_X)$ is a flat $n$-dimensional torus, $B\in H^2(X,{\R}/{\Z})$
and $\phi_X^{triv}$ is the constant map form $X$ to the trivial theory point in the moduli space of CFTs.
These data in turn give rise to a toroidal CFT, which can be realized
as a sigma model with the target space $(X,g_X)$ and
given B-field $B$. The residual bundle of CFTs on $X$ is trivial. 

Let us consider a $1$-parameter family of CFTs defined by the family
$(X,\lambda g_X, \phi_X^{triv};B=0)$, where $\lambda \in (0,+\infty)$.
There are two degenerations of this family, which define
two points of the boundary $\partial \overline{\cal M}_{c=n}^{tor}$.
 As $\lambda \to +\infty$,
we get a toroidal CFT defined by $(X,{\R}_+^{\ast}\cdot g_X, \phi_X^{triv};B=0)$.
As $\lambda \to 0$ we get $(X^{\vee},{\R}_+^{\ast}\cdot g_{X^{\vee}}, \phi_X^{triv};B=0)$,
where $(X^{\vee},g_{X^{\vee}})$ is the dual flat torus.

There might be further degenerations of the lattice.
Thus one obtains a stratification of the compactified moduli space
of lattices (and hence  CFTs).
Points of the compactification are described by flags of vector spaces
$0=V_0\subset V_1\subset V_2\subset...\subset V_k\subset {\R}^n$.
In addition one has a lattice $\Gamma_{i+1}\subset V_{i+1}/V_i$,
considered up to a scalar factor. These data give rise to
a tower of torus bundles $X_k\to X_{k-1}\to...\to X_1\to \{pt\}$
over tori with fibers $(V_{i+1}/V_i)/\Gamma_{i+1}$.
If $V_k\simeq {\R}^{n-l}, l\ge 1$, then one has also
a map from the total space $X_k$
of the last torus bundle to the point $[H_k]$ in the moduli space
of toroidal theories of smaller central charge:
 $\phi_n:X_k\to {\cal M}_{c=l}^{tor}$,
$\phi_k(X_k)=[H_k]$.

\subsection{Example: WZW model for $SU(2)$}
 
In this case we have a discrete family with
 $c={3k\over{k+2}}$, where $k\ge 1$ is an integer number
called {\it level}.
In the limit $k\to +\infty$ one gets $X=SU(2)=S^3$
equipped with the standard metric.
The corresponding bundle is the trivial bundle of
trivial CFTs (with $c=0$ and
$H=H^{0,0}={\C}$). Analogous picture holds for an arbitrary
compact simply connected simple group $G$.

\subsection{A-model and B-model of $N=2$ SCFT as boundary strata}

The boundary of the compactified moduli space
$\overline {{\cal M}}^{N=2}$ of $N=2$ SCFTs with a given central charge
contains an open
stratum given by sigma models with Calabi-Yau targets.
Each stratum is parametrized by the classes of equivalence of
quadruples $(X,J_X,{\R}_+^{\ast}\cdot g_X,B)$ where $X$ is
a compact real manifold, $J_X$ a complex structure,
$g_X$ is a Calabi-Yau metric, 
 and $B\in H^2(X,i{\bf R}/{\bf Z})$ is a
$B$-field. The residual bundle of CFTs is a bundle of free fermion
theories.

 As a consequence of supersymmetry,
the moduli 
space ${\cal M}^{N=2}$  of superconformal field theories is a complex
manifold which 
is locally
isomorphic to the product of two complex manifolds. \footnote{Strictly speaking,
one should exclude models with chiral fields of conformal dimension $(2,0)$,
 e.g. sigma models on hyperk\"ahler manifolds, see [AM].}
It is believed that this decomposition (up to certain
corrections) is global.  
Also, there are two types of sigma models with Calabi-Yau targets:
$A$-models and $B$-models.
Hence, the traditional picture
of the compactified moduli space looks as follows:

\vspace{3mm}

\centerline{\epsfbox{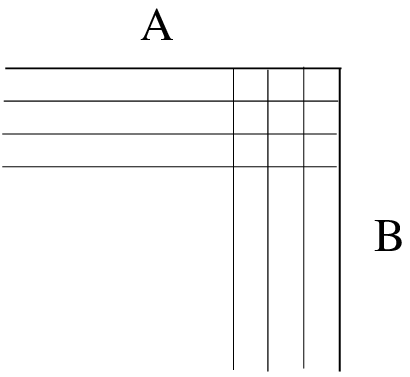}}

\vspace{2mm}

Here the boundary consists of two open strata (A-stratum and B-stratum)
and a mysterious meeting point. This point corresponds,
in general, to a submanifold
of codimension one in the closure of A-stratum and of B-stratum.

We argue that this picture should be modified. There is another open
stratum of $\partial \overline {{\cal M}}^{N=2}$ (we call it T-stratum).
 It consists of toroidal
 models (i.e. CFTs  associated with Narain lattices), 
 parametrized by a 
 manifold $Y$ with a Riemannian metric defined up to a scalar factor.
 This subvariety meets both $A$ and $B$ strata along the codimension
 one stratum corresponding to the double collapse. 
 Therefore the ``true'' picture
 is obtained from the traditional one
 by the real blow-up at the corner:
 
 \vspace{3mm}

\centerline{\epsfbox{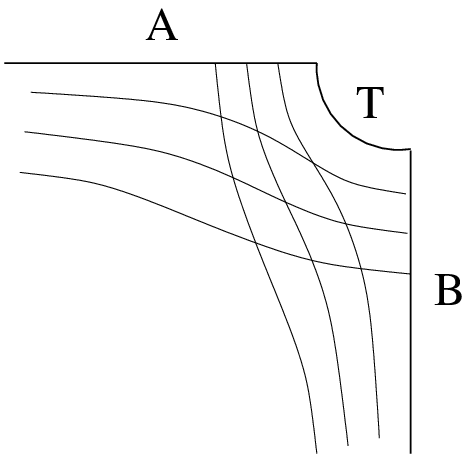}} 
 
\subsection{Mirror symmetry and the collapse}

Mirror symmetry is related to the existence
of two different strata of the boundary $\partial \overline {{\cal M}}^{N=2}$
which we called A-stratum and B-stratum. 
As a corollary, same quantities
admit different geometric descriptions near different
strata. In the traditional picture, one can introduce
natural coordinates in a small neighborhood of a boundary point
corresponding to $(X,J_X,{\R}_+^{\ast}\cdot g_X,B)$.
Skipping $X$ from the notation, one can say that the coordinates
are $(J,g,B)$ (complex structure, Calabi-Yau metric and the B-field).
Geometrically, the pairs $(g,B)$ belong to the preimage of the
K\"ahler cone under the natural map $Re: H^2(X,{\C})\to H^2(X,{\R})$
(more precisely, one should consider $B$ as an element
of $H^2(X,i{\R}/{\Z})$). It is usually said, that one considers
an open domain in the complexified K\"ahler cone with the property that
with the class of metric $[g]$ it contains also the ray
$t[g], t\gg 1$. The mirror symmetry gives rise to an
identification of neighborhoods of $(X,J_X,{\R}_+^{\ast}\cdot g_X,B_X)$
and $(X^{\vee},J_{X^{\vee}},{\R}_+^{\ast}\cdot g_{X^{\vee}},B_{X^{\vee}})$
such that $J_X$ is interchanged with $[g_{X^{\vee}}]+iB_{X^{\vee}})$
and vice versa. 

We can describe this picture in a different
way.
Using the identification of complex and K\"ahler moduli, one can choose 
$([g_X],B_X, [g_{X^{\vee}}],B_{X^{\vee}})$ as local coordinates near
the meeting point of A-stratum and B-stratum.
 There is an action of the additive
semigroup ${\R}_{\ge 0}\times {\R}_{\ge 0}$
in this neighborhood. It is given explicitly by the formula
$([g_X],B_X, [g_{X^{\vee}}],B_{X^{\vee}})\mapsto 
 (e^{t_1}[g_X],B_X, e^{t_2}[g_{X^{\vee}}],B_{X^{\vee}})$ where
 $(t_1,t_2)\in {\R}_{\ge 0}\times {\R}_{\ge 0}$.
 As $t_1\to +\infty$, a point of the moduli space
 approaches  the B-stratum,
 where the metric is defined up to a positive scalar only.
 The action of the second semigroup ${\R}_{\ge 0}$ extends
by continuity to the non-trivial action on the B-stratum.
Similarly, in the limit
 $t_2\to +\infty$ the flow retracts the point to the A-stratum.

This picture should be modified, if one makes a real blow-up at
the corner, as discussed before. Again, the action of the semigroup
${\R}_{\ge 0}\times {\R}_{\ge 0}$ extends continuously to the
boundary.
Contractions to A-stratum and B-stratum carry non-trivial
actions of the corresponding semigroups isomorphic to 
${\R}_{\ge 0}$. Now, let us choose a point
in, say, A-stratum. Then the semigroup flow takes it 
along the boundary to
the new stratum, corresponding to the double collapse.
The semigroup ${\R}_{\ge 0}\times {\R}_{\ge 0}$ acts trivially on this stratum.
A point of the double collapse is also
a limiting point of a $1$-dimensional orbit of
${\R}_{\ge 0}\times {\R}_{\ge 0}$ acting on the T-stratum.
Explicitly, the element $(t_1,t_2)$ changes the size 
of the tori defined by the Narain lattices, rescaling them with the coefficient
$e^{t_1-t_2}$. This flow carries the point of T-stratum 
 to another point
of the double collapse, which can be moved then inside of the B-stratum.
The whole path, which is the intersection of $\partial \overline{\cal M}^{N=2}$
and the ${\R}_{\ge 0}\times {\R}_{\ge 0}$-orbit,
connects an A-model with the corresponding B-model
through the stratum of toroidal models. We can depict it as
follows:

\vspace{3mm}

\centerline{\epsfbox{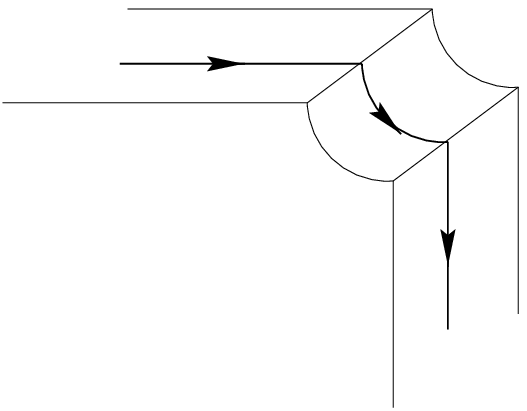}}

\vspace{2mm}

The T-portion of the path (we call it T-path) connects 
dual torus fibrations over the same Riemannian base.
This is mirror symmetry in our picture.

This description is inspired by [SYZ]. The reader 
notices however,
that in our picture, the mirror symmetry phenomenon is explained entirely 
in terms of the {\it boundary} of the compactified moduli space.
In order to explain the mirror symmetry phenomenon
it is not necessary to build full SCFTs. 
It is sufficient to work with simple
toroidal models on the boundary of the compactified moduli
space $\overline {{\cal M}}^{N=2}$. Also, in contrast with [SYZ], we do not
use supersymmetric cycles (D-branes) in our description.

\section{Calabi-Yau manifolds in the large complex structure limit}

\subsection{Maximal degenerations of Calabi-Yau manifolds}

Let ${\C}_q^{mer}=\{f=\sum_{n\ge n_0}a_nq^n\}$ be the field of germs
at $q=0$ of meromorphic functions in one complex variable.

Let ${\cal X}_{mer}$ be an algebraic $n$-dimensional Calabi-Yau manifold 
over ${\C}_q^{mer}$ (i.e. ${\cal X}_{mer}$ is a smooth projective manifold
over ${\C}_q^{mer}$ with the trivial canonical class: $K_{{\cal X}}=0$).
We fix an algebraic non-vanishing
volume element $vol\in \Gamma ({\cal X}_{mer},K_{{\cal X}})$.
The pair $({\cal X}_{mer},vol)$ defines a 1-parameter 
analytic family of complex
Calabi-Yau manifolds $(X_q,vol_q), 0<|q|<r_0$, for  some $r_0>0$.

Let $[\omega]\in H^2_{DR}({\cal X}_{mer})$ be the cohomology class in the ample cone.
Then for every $q$, such that $0<|q|<r_0$ it defines a K\"ahler class
$\omega_q$ on $X_q$. By the Yau theorem, there exists a unique Calabi-Yau
metric $g_{X_q}$ on $X_q$ with the K\"ahler class $[\omega_q]$.

 It follows from the resolution of singularities,
that as $q\to 0$ one has the following formula:

$$\int_{X_q}vol_q\wedge \overline{vol}_{q}=C(log|q|)^m|q|^{2k}(1+o(1))$$
for some $C\in {\C}^{\ast}, k\in {\Z}, 0\le m\le n=dim\,({\cal X}_{mer})$.

\begin{dfn} We say that ${\cal X}_{mer}$ has maximal degeneration at $q=0$
if in the formula above we have $m=n$.

\end{dfn}

Let us show that this definition is equivalent to the usual one,
given in terms of variations of Hodge structures (see [Mo]).

\begin{prp} The Calabi-Yau manifold
${\cal X}_{mer}$ has maximal degeneration iff 
for all sufficiently small $q$ there exists
          a vector $v\in H^n(X_q,\C)$ such that $(T-Id)^{n+1}v=0$ and
           $(T-id)^n v\ne 0$ where $T$ is the
         monodromy operator.

\end{prp}

 {\it Proof.}
 First of all, notice that the volume
  $\int_{X_q}vol_q\wedge \overline{vol}_q$ can be calculated cohomologically
      as the Poincar\'e pairing $\langle [vol_q],  [\overline{vol}_q]\rangle$
in the (primitive part of) middle cohomology $H^n(X_q,\C)$.

 We can assume (after passing to a cover by adding a root of $q$) that
 the operator
  $T$ is unipotent. Let us trivialize the bundle $H^n(X_q,\C)$ over the
  punctured disc
   by multplication by
   $$q^{-log(T)/2\pi i}=\sum_{k=0}^n \left({log(q)\over -2\pi i}\right)^k
   {(log(T))^k\over k!}$$
 The nilpotent orbit
   theorem says that the Hodge filtration extends to a holomorphic filtration
    on the trivialized bundle over the whole disc, including the point $q=0$.
    Thus, the bundle $q^{-log(T)/2\pi i}(H^{n,0})$ extends to a line bundle over
the disc,
     and the section  $q^{-log(T)/2\pi i}([vol_q])$  is a non-zero
      meromorphic section of
     this bundle.
     After the multiplication by an appropriate power of $q$ we may assume
further
     that this section is holomorphic and non-vanishing at $q=0$.  We denote
     this holomorphic section by $a(q),\,\,\,a(0)\ne 0$.

     Now let us calculate the volume:
     $$\langle [vol_q],  [\overline{vol}_q]\rangle=
     \sum_{k,l\ge 0}\left\langle \left({log(q)\over 2\pi i}\right)^k
     {(log(T))^k\over k!}a(q), \left({log(\overline{q})\over -2\pi i}\right)^l
     {(log(T))^l\over l!}\overline{a(q)}\right\rangle=$$
     $$= \sum_{k,l\ge 0} \left({log(q)\over 2\pi i}\right)^k
     \left({log(\overline{q})\over -2\pi i}\right)^l {(-1)^l\over k!l!}
      \left\langle (log(T))^{k+l} a(q),\overline{a(q)}\right\rangle$$
     Here we use the fact that operator $log(T)$  is real and also
skew-symmetric
      with respect to the Poincar\'e pairing.
       It follows from the  equality $log(T)^{n+1}=0$ (which holds
automatically)
       that in the sum above all
        terms with $k+l> n$ vanish. The contribution of
terms with $k+l=n$
         is equal to
         $${(log(|q|)^n\over (\pi i)^n n!}\langle (log(T)^n a(q),
         \overline{a(q)}\rangle.$$  
          It follows easily from the standard
properties
          of variations of polarized Hodge structures (see e.g. [De]) that
           if $log(T)^n\ne 0$ then  in fact $log(T)^n(a(0))\ne 0$,
            and moreover $\langle (log(T)^n a(0),
         \overline{a(0)}\rangle\ne 0$.
         Now the statement of the Proposition becomes obvious.
$\blacksquare$

Notice that in [De] a slightly stronger condition 
of maximal degeneration was  imposed:
 the weight filtration on $H^*(X_q,\C)$ associated with the monodromy operator
  should be complementary to the Hodge filtration.

Let us recall the definition of the Gromov-Hausdorff metric $\rho_{GH}$.
It is a metric on the space of isometry classes of compact
metric spaces. We say that two metric spaces $M_1$ and $M_2$ are
$\varepsilon$-close in $\rho_{GH}$ if there exists a metric space
$M$ containing both $M_1$ and $M_2$ as metric subspaces, such that
$M_1$ belongs to the $\varepsilon$-neighborhood of $M_2$ and
vice versa. Then $\rho_{GH}(M_1,M_2)$ is given by the minimum
of such $\varepsilon$.

Let us rescale the Calabi-Yau metric: 
$g_{X_q}^{new}=g_{X_q}/diam(X_q,g_{X_q})^{2}$.
Thus we obtain a 1-parameter family of Riemannian manifolds
$X_q^{new}=(X_q,g_{X_q}^{new})$ of the diameter $1$.

\begin{conj} If ${\cal X}_{mer}$ has maximal degeneration at $q=0$ then
there is a limit $(\overline{Y},g_{\overline{Y}})$ of $X_q^{new}$
in the Gromov-Hausdorff metric, such that:

a) $(\overline{Y},g_{\overline{Y}})$ is a compact metric space,
which contains a smooth oriented 
Riemannian manifold $(Y,g_Y)$ of dimension $n$
as a dense open metric subspace.
The Hausdorff dimension of $Y^{sing}=\overline{Y}\setminus Y$
is less or equal than $n-2$.

b) $Y$ carries an {\it integral affine structure}.
This means that it carries a torsion-free
flat connection $\nabla$ with the holonomy 
contained in $SL(n,{\Z})$.

c) The metric $g_Y$ has a potential. This means that it is locally given 
in affine coordinates by a symmetric matrix
$(g_{ij})=(\partial^2 K/\partial x_i\partial x_j)$, where $K$ is a
smooth function (defined modulo adding an affine function, i.e. the sum of a linear function and a constant).

d) In affine coordinates the metric volume element is constant,
 $det(g_{ij})=det(\partial^2 K/\partial x_i\partial x_j)=const$ 
 (real Monge-Amp\`ere equation).

\end{conj}

At the end of this section we propose a non-rigorous explanation
of our conjecture based on differential-geometric considerations.

\begin{rmk} 

1) Since the matrix $(g_{ij})$ defined by the
metric $g_Y$ is positive, the function $K$ is convex. 
In particular, there is locally well-defined Legendre transform
of $K$. This fact will be used later, when we will discuss the duality
of Monge-Amp\`ere manifolds.

2) It seems plausible that in the case when all
$X_q$ are simply-connected, and $h^{k,0}(X_q)=0$ for $0<k<n$,
the metric space $\overline{Y}$ is a homological
sphere of dimension $n$. In all examples it is
in fact homeomorphic to $S^n$.

\end{rmk}

The conjecture opens
the way for compactification of
 the moduli space of Calabi-Yau metrics on a given Calabi-Yau
manifold $M$, by adding as a boundary component the set of pairs 
$(\overline{Y},{\R}_+^{\ast}\cdot g_{\overline{Y}})$ for all
 1-parameter maximal degenerations $X_q$, such that $X_{q^{\prime}}=M$
 for some $q^{\prime}$.
This corresponds to a choice of
 a ``cusp'' in the moduli space of Calabi-Yau manifolds.
This choice is usually described
in terms of certain algebro-geometric
data: the action of the monodromy operator, 
variation of Hodge structures, mixed Hodge
structure of the special fiber, etc.
The previous conjecture
offers a pure ``metric'' description of a cusp.
 
It follows from part b) of the conjecture that one can choose a 
$\nabla$-covariant lattice $T_{Y,y}^{\Z}\subset T_{Y,y}, y\in Y$.
Suppose we are given a triple $(Y,g_Y,\nabla)$, satisfying the 
properties a)-c) of the conjecture, and we have fixed
a covariant lattice $T_{Y}^{\Z}$ in
the tangent bundle $T_Y$. Then we can construct a
1-parameter family of {\it non-compact} complex Calabi-Yau manifolds,
endowed with Ricci flat K\"ahler metrics.
Namely, let $X^{\varepsilon}$ be the total space of the torus bundle
$p_{\varepsilon}:X^{\varepsilon}\to Y$ with fibers
$T_{Y,y}/\varepsilon T_{Y,y}^{\Z},
y\in Y, 0<\varepsilon\le \varepsilon_0$. The total space
$TY$ of the tangent bundle $T_Y$ carries a canonical complex
structure coming from the isomorphism 
$T_{TY}\simeq \pi^{\ast}T_Y\oplus \pi^{\ast}T_Y
\simeq \pi^{\ast}T_Y\otimes {\C}$ where $\pi:TY\ra Y$ is
the canonical projection (here we use the affine structure on $Y$).
Using the same identification, we introduce a metric on $TY$, namely
$g_{TY}=\pi^{\ast}g_Y\oplus \pi^{\ast}g_Y$. It is easy to see, that
$g_{TY}$ is a K\"ahler metric with the potential $\pi^{\ast}K$. It follows from the
Monge-Amp\`ere equation that the metric $g_{TY}$ is Ricci flat.
Passing
to the quotient, we obtain on $X^{\varepsilon}$
 a complex structure $J_{X^{\varepsilon}}$
and a Ricci flat K\"ahler metric $g_{X^{\varepsilon}}$.

Let $U\subset Y$ be an open simply-connected subset.
Then there is an action of the torus $T^n\simeq T_{Y,y}/T_{Y,y}^{\Z}$
on $p_{\varepsilon}^{-1}(U), y\in U$ (different tori are
identified for different points $y\in U$ by means of the 
connection $\nabla$). It implies that for any 
$t\in H^1(Y,(T_Y/T_Y^{\Z})^{discr})$ (cohomology with coefficients in
the local system of tori considered as abstract groups) one
can define a twisted manifold $X^{\varepsilon,t}$, which is
the total space of the torus fibration 
$p_{\varepsilon,t}: X^{\varepsilon,t}\to Y$.

Roughly speaking, the next conjecture
says that the ``leading asymptotic term'' of the family of Calabi-Yau
 manifolds $X_q^{new}, q=e^{-1/\varepsilon}$
 near the point of maximal degeneration $\varepsilon=0$, 
is isomorphic up to a twist to the family
 $(X^{\varepsilon},J_{X^{\varepsilon}})$ associated
with the torus bundle described above.

More precisely, we formulate it as follows.

\begin{conj} Let $({\cal X}_{mer}, vol)=(X_q, vol_q)$ 
be a 1-parameter family of maximally degenerate
Calabi-Yau manifolds, and $X_q^{new}$ be the family with rescaled metrics,
as before. There exist a constant $C>0$ and a function $t(q)$ such that
K\"ahler manifolds  $X_q^{new}$ and 
$X^{\varepsilon (q),t(q)}$ with $\varepsilon (q)=C(log|q|)^{-1}$ are close
to each other (as $q\to 0$) in the following sense:

for any $\delta>0$ there exist a decomposition 
$X_q=X_q^{sm}\sqcup X_q^{sing}$ and an embedding of smooth manifolds
$j_q:X_q\to p_{\varepsilon(q),t(q)}^{-1}(Y\setminus (Y^{sing})^{\delta})$, where
$(Y^{sing})^{\delta}$ is a $\delta$-neighborhood of $Y^{sing}$, such
that:

a) $(X_q,X_q^{sing})$ converges in the Gromov-Hausdorff metric
to the pair $(\overline{Y},Y^{sing})$.

b) $j_q$ identifies up to $o(1)$ terms, uniformly in $x\in X_q^{sm}$,
the scalar products and complex structures on 
the tangent spaces $T_x X_q$ and
$T_{j_q(x)}X^{\varepsilon (q),t(q)}$.

\end{conj}

There is the following motivation for the Conjectures 1 and 2.
In general, for a degenerating family of Riemannian metrics 
with non-negative Ricci
curvature, one expects a description in terms of a
tower of fibrations (collapses) with singularities (compare with 2.3).
 \footnote { Some 
 steps in the program of compactification of the space of metrics
  are accomplished now (see e.g. [CC]), 
  but still there are many non-clarified issues.}
  In the case of K\"ahler manifolds there are two basic pictures of a simple
collapse. The first case is when both the base and the fiber are 
K\"ahler manifolds. In the second case fibers are flat totally real tori
of dimension $m$ and the base looks locally as a product of a domain
in ${\R}^m$ with a K\"ahler manifold.
The logarithmic factor in the asymptotic behavior of the volume should come
only from torus fibers. Thus, the largest possible power of the
logarithm can appear only when we have a tower of purely torus fibrations.
It seems that the fixing (up to a scalar) of  the K\"ahler class forbids the
multiple collapse. These considerations give an intuitive
``explanation'' of our conjectures.

\begin{rmk} During the preparation of this text we learned
that  conjectures similar to ours were proposed
 independently by M.~Gross and P.~Wilson (see [GW]). A remarkable
  achievement in
  [GW] consists of the verification of conjectures
  in the case
    of degenerating $K3$ surfaces, together with a precise
    description of the behavior of metrics near singular fibers.
  Also, in a recent preprint [Le]
 mirror symmetry was discussed from a similar point of view. 
 In the main body of the present paper
 we will consider degenerations of complex abelian varieties.
  In this case the conjectures obviously hold.

\end{rmk}

\subsection{Monge-Amp\`ere manifolds and duality of torus fibrations}

In this section we propose a mathematical language for the geometric mirror symmetry,
understood as a duality of torus fibrations.

\begin{dfn} A Monge-Amp\`ere manifold is a triple $(Y,g,\nabla)$,
where $(Y,g)$ is a smooth Riemannian manifold with the metric $g$, and
$\nabla$ is
a flat connection on $T_Y$ such that:

a) $\nabla$ defines an affine structure on $Y$.

b) Locally in affine coordinates $(x_1,...,x_n)$
 the matrix $((g_{ij}))$
 of $g$ is given by $((g_{ij}))=((\partial^2 K/\partial x_i\partial x_j))$
 for some smooth real-valued function $K$.
 
c) The Monge-Amp\`ere equation
$det((\partial^2 K/\partial x_i\partial x_j))=const$ is satisfied.

\end{dfn}

Monge-Amp\`ere manifolds were studied (under a different name)
in [CY] where is was proven that
 if $Y$ is {\it compact} then its finite cover is  a torus.

Let us consider a (non-compact) example motivated by the mirror symmetry for K3 surfaces
 (see also [GW]).
Let $S$ be a complex surface endowed with a holomorphic
non-vanishing volume form $vol_S$, and $\pi:S\to C$ be a holomorphic
fibration over a complex curve $C$, such that fibers of $\pi$
are non-singular elliptic curves.

We define a metric $g_C$ on $C$ as the K\"ahler metric associated
with the $(1,1)$-form $\pi_{\ast}(vol_S\wedge \overline{vol}_S)$.
Let us choose (locally on $C$) a basis $(\gamma_1,\gamma_2)$
in $H_1(\pi^{-1}(x),{\Z}), x\in C$. We define two closed 1-forms on $C$ by the
formulas

$$\alpha_i=Re(\int_{\gamma_i}vol_S), i=1,2.$$

It follows that $\alpha_i=dx_i$ for some
functions $x_i, i=1,2$. We define an affine structure on $C$,
and the corresponding connection $\nabla$,
by saying that $(x_1,x_2)$ are affine coordinates.
One can check directly that $(C,g_C,\nabla)$ is a Monge-Amp\`ere
manifold. In a typical example of elliptic fibration
of a K3 surface, one gets $C={\C}{\bf P}^1\setminus \{z_1,...,z_{24}\}$,
where $\{z_1,...,z_{24}\}$ is a set of distinct $24$ points in 
${\C}{\bf P}^1$.

Returning to the general case, we can restate a portion 
of our conjectures by saying
that the smooth part of the Gromov-Hausdorff
limit of a maximally degenerate family of Calabi-Yau manifolds
is a Monge-Amp\`ere manifold with an integral affine structure.

There is a well-known duality on local solutions of the Monge-Amp\`ere
equation.

\begin{lmm} Let $U\subset {\R}^n$ be a convex open domain 
in ${\R}^n$ equipped with the standard affine coordinates $(x_1,...,x_n)$,
and $K:U\to {\R}$ be a convex function satisfying the Monge-Amp\`ere
equation. Then the Legendre transform $\widehat{K}(y_1,...,y_n)=
max_{x\in U}(\sum_ix_iy_i-K(x_1,...,x_n))$ also satisfies the
Monge-Amp\`ere equation.

\end{lmm}

{\it Proof.} The graph of $L=dK$ is a Lagrangian
submanifold in $T^{\ast}{\R}^n={\R}^n\oplus ({\R}^n)^{\ast}$.
Let $p_1$ and $p_2$ be the natural projections to the direct
summands. They are local diffeomorphisms.
Since $K$ is defined up to the adding of an affine function,
the graph itself is defined up to translations.
The Monge-Amp\`ere equation corresponds to the condition
$p_1^{\ast}(vol_{{\bf R}^n})=p_2^{\ast}(vol_{{\bf R}^{\ast n}})$,
where $vol_{{\bf R}^n}$ (resp. $vol_{{\bf R}^{\ast n}}$)
denotes the standard volume form on ${\bf R}^n$ (resp. ${\bf R}^{\ast n}$). 
The manifold $L$ can be considered as a graph of $d\widehat{K}$.
Thus $\widehat{K}$ satisfies the Monge-Amp\`ere equation as well.
The Lemma is proved. $\blacksquare$

The manifold $L$ carries a Riemannian metric $g_L$ induced
by the indefinite metric $\sum_idx_i\,dy_i$ on ${\R}^n\oplus ({\R}^n)^{\ast}$.
This metric is given by the matrix $(\partial^2 K/\partial x_i\partial x_j)$
in coordinates $(x_1,...,x_n)$, and by the matrix
$(\partial^2 \widehat{K}/\partial y_i\partial y_j)$ in the dual coordinates.
Thus on $L$ we have a metric, and two affine structures
(pullbacks of the standard affine structures on the coordinate spaces).
Hence we have two structures of the Monge-Amp\`ere manifold on $L$.
It is easy to see that the local pictures can be glued together.
This leads to the following result.

\begin{prp} For a given Monge-Amp\`ere manifold $(Y,g_Y,\nabla_Y)$
there is a canonically defined
dual Monge-Amp\`ere manifold $(Y^{\vee},g_Y^{\vee},\nabla_Y^{\vee})$
such that $(Y,g_Y)$ is identified
with $(Y^{\vee},g_Y^{\vee})$ as Riemannian manifolds, and
the local system $(T_{Y^{\vee}},\nabla_Y^{\vee})$ 
is naturally isomorphic to the local system dual to  
$(T_Y,\nabla_Y)$.

\end{prp}

\begin{cor} If $\nabla_Y$ defines an integral affine structure on $Y$
(i.e. the holonomy of $\nabla_Y$ belongs to $GL(n,{\Z})$),
then $\nabla_Y^{\vee}$ defines an integral affine structure on $Y^{\vee}$.
As the dual covariant lattice one takes the lattice $(T^{\Z}_Y)^{\vee}$,
which is dual to $T_Y^{\Z}$ with respect to the metric $g_Y$.

\end{cor}

Now we can state the geometric counterpart of the mirror symmetry conjecture.

\begin{conj} Smooth parts of maximal degenerations of dual families
of Calabi-Yau manifolds are dual Monge-Amp\`ere manifolds
with dual integral affine structures.

\end{conj}

Monge-Amp\`ere manifolds with integral affine structures are real
analogs of Calabi-Yau manifolds. In fact the mirror duality
in the sense of this section holds for a larger class of
manifolds.
We define an {\it AK-manifold} (AK stands for {\it affine} and {\it K\"ahler})
as in the Definition 2, but dropping the condition c) (Monge-Amp\`ere equation),
 see also [CY].
The reader can check easily that all constructions of this section,
including the duality of torus fibrations hold for AK-manifolds as well.

\begin{rmk} The idea to use the Legendre transform for the
purposes of mirror symmetry was around for some time
(see for example [H], [Le]).

\end{rmk}

\begin{rmk} In our description of geometric mirror symmetry we
ignore the B-fields. In what follows we will always assume
that $B=0$.

\end{rmk}

\subsection{Speculations about relations with non-archimedean geometry}

Considerations from CFT and from differential geometry
indicate that the  integral affine structure on $Y$ 
does not depend on the choice of the K\"ahler class of Calabi-Yau metrics.
Thus, we obtain a ``combinatorial'' invariant $(Y,T^{\Z}_Y)$
of (maximally degenerating)
 Calabi-Yau variety over the local field $K={\C}((q))$.      
One can argue that in this case
  there will be a canonical atlas of coordinate charts such that the transition maps belong
   to the group $SAff(n,\Z):=SL(n ,\Z)\ltimes \Z^n$.
  The natural question arises whether one can define and calculate it
   purely algebraically, without the use of transcendental methods
    and Calabi-Yau metrics.  We expect that the answer to this question
    is positive. In other words
     there exists a canonical way to associate the data 
    $(Y,T^{\Z}_Y)$
   with arbitrary smooth projective variety $X,\,c_1(T_X)=0$ 
   having ``maximal degeneration'' over an arbitrary
     field $K$  with  a discrete valuation.
     
     The conjectural answer 
     (only for the compactification $\overline{Y}$ of $Y$)
     is the following: let us choose
      (after an extension of the field $K$) a model with stable reduction.
      Call an irreducible component $D$ of the special fiber $X_0$
      {\it  essential} if the order of pole at $D$ of the global volume element on $X$
      is maximal among all components of $X_0$.
       We define topological space
       $\overline{Y}(X_0)$ as the Clemens complex spanned by essential divisors
        (see [LTY]).
       Roughly speaking, $k$-cells of $\overline{Y}$ correspond 
       to irreducible components
        of $(k+1)$-fold intersections of essential divisors.
        
         Using ideas from motivic integration and
        from Berkovich theory of non-archimedean analytic spaces (see [Be]),
        one can prove the following result.         
        
        \begin{prp}For
          different choices of models with stable reduction, the spaces 
          $\overline{Y}(X_0)$ can be canonically identified.
        
        \end{prp}
        
 We will sketch the idea of the proof.
 \footnote{After the lecture of M.K. at Rutgers
 University on October 30, 2000.}
Following [Be], one can associate to an algebraic variety $X$ over
a local field $K$,
a $K$-analytic space $X^{an}$. For an affine $X$, points of $X^{an}$
 are ${\R}$-valued multiplicative
seminorms on the algebra of functions ${\cal
O}(X)$, extending a given norm on $K$. In our case $K={\C}((q))$. 
 
{\it Step 1}. One defines the set of {\it divisorial points} $D(X,K)$
as the direct limit of the sets of irreducible components of special
fibers, taken over all models of $X$ over ${\C}[[q]]$ (not necessarily
smooth). One has a natural
embedding $D(X,K)\subset X^{an}$ (the 
seminorm corresponding to $D$ is the exponent of the order of the pole of the 
volume form on the divisor $D$).  

Similarly one defines 
$D(X,K^{\prime})=D(X\otimes K^{\prime})/Gal(K^{\prime}/K)$,
where $K^{\prime}/K$ is a Galois extension of finite order.
Taking direct limit of $D(X,K^{\prime})$ over all $K^{\prime}$
one gets the set $D(X)$.

{\it Step 2}. One checks that $D(X)\subset X^{an}$, and in fact $D(X)$
is a dense subset in $X^{an}$. For local fields with finite or countable
 residue field the set $D(X)$
is  countable, while $X^{an}$ is connected.

{\it Step 3}. Having a divisor $D\in D(X)$ one defines the number
$$v(vol,D)={ord_{(X\otimes K^{\prime},D)}(vol)+1\over{ram(K^{\prime}/K)}},$$
where $ord$ means the order of zero of the volume form $vol$
on $D\subset X\otimes K^{\prime}$, and $ram$ means
the ramification index. One checks that $v(vol,D)$ is well-defined.
In $p$-adic case this number can be derived from the contribution
 of $D$  to the integral 
$\int_{X(K)}|vol|_K$, where $|vol|_K$ is the associated measure.

{\it Step 4}. One checks that $\overline{Y}(X_0)$ is the closure
in $X^{an}$ of the set of divisors $D\in D(X)$ such that $v(vol,D)$
is minimal. Hence it is independent of a choice of a model.
This explains the statement of the theorem. $\blacksquare$
         
           In examples coming from toric geometry the space 
           $\overline{Y}=\overline{Y}(X_0)$ is always a manifold.
           It is not clear yet 
           what is the origin of the smooth 
           part $Y\subset \overline{Y}$, and of the affine structure
            on it.    Conjecturally, all this
             comes from a map 
             $\pi:X(\overline{K})\to   \overline{Y}$ where $\overline{K}$ 
             is the
            algebraic closure of $K$. 
            In the differential-geometric picture 
            of torus fibrations (when $K=\C_{mer}$)
             the map   $\pi$ is obvious:
             it associates with a meromorphic (finitely ramified)  
             family of points $x_q\in X_q$
             the limit point $lim_{q\to 0}\, x_q\in Y$  in the metric sense. 
             Also, the differential-geometric picture suggests
              that the closure of the image $\pi(Z(\overline{K}))$ 
              where $Z\subset K$ is an algebraic  subvariety, should
              be a piecewise linear closed subset of $Y$, 
              and linear pieces of it have rational directions.
              In particular, if $Z$ is a curve then  $\pi(Z(\overline{K}))$
               is a graph in $Y$. This opens a way to express Gromov-Witten 
                 invariants of $X$ in terms of the
               Feynman expansion for certain quantum field theory on $Y$.
           
           Also, we expect that the choice of 
           an ample class in $NS(X)\otimes\R$
           on $X$ gives rise to the {\it dual }
            integral affine structure on $Y$
        defined again in some purely algebro-geometric way. 
        If the ample
         class is the first Chern class
         of a line  bundle, then 
         there should be also a canonical reduction of the dual integral
          affine structure to a $SAff(n,\Z)$ -structure.

\section{$\A$-algebras and $\A$-categories}

\subsection{Two problems with the general definition}

The purpose of this section is to describe the framework in which
the results concerning $\A$-categories
will be formulated. We would like to make few comments 
even before recalling a definition of the Fukaya category.
These comments are informal. Precise definitions will be given
later in this section.

There are two main problems with the definition
of the Fukaya category. First, morphisms
can be defined only for transversal Lagrangian submanifolds
 (in particular, the identity morphism is never defined). Second,
since there are pseudo-holomorphic discs with the boundary on a given
Lagrangian submanifold, one has to add a composition $m_0$ to the
set of compositions $m_n, n\ge 1$. As a result, the spaces of
morphisms are not complexes: $m_1^2\ne 0$. On the other hand, the
derived category of coherent sheaves arises from an
$\A$-category without $m_0$ and with
the condition $m_1^2=0$. 
Hence one should explain in which sense two $\A$-categories in question
are equivalent. 
The above-mentioned problems can be resolved by an appropriate
generalization of the notion of $\A$-category. We offer such a generalization
below (we call it $\A$-pre-category).
We remark that there exists a better generalization.
It involves numerous preparations and will be given elsewhere (see [KoS]).
On the other hand, the problem with $m_0$ does not
appear in the version of the Fukaya
category for abelian varieties considered in this paper.
Hence, for the purposes
of present paper it is sufficient to work with $\A$-pre-categories
(or $\A$-categories with transversal structure, cf. [P1]).
This gives a partial solution to the transversality problem,
and provides a solution to the problem with the identity
morphisms.

Using $\A$-pre-categories we are going to formulate and prove in Section 8
a variant of the Homological Mirror Conjecture
for torus fibrations.
It can be applied to the case of abelian varieties.
In particular, one can obtain certain formulas for Massey
products for abelian varieties in terms of partial
theta-sums similar to those considered in [P1].

\subsection{Non-unital $\A$-algebras and $\A$-categories}

Let $A=\oplus_{i\in \Z}A^i$ be a
 $\Z$-graded module over a field $k$. \footnote{In what follows one
can replace $k$  by a $\Z$-graded commutative associative algebra and assume
 that all $k$-modules are projective.}
As usual, we will denote by $A[n]$ the graded $k$-module such that
$(A[n])^i=A^{i+n}$ for all $i$.

\begin{dfn}  A structure of non-unital $\A$-algebra on $A$ is given
by a codifferential $d$ of degree $+1$ on the cofree tensor coalgebra
$T_+(A[1])=\oplus_{n\ge 1}(A[1])^{\otimes n}$.

\end{dfn}

The codifferential $d$ is by definition a coderivation, such that $d^2=0$.
It is uniquely determined by its ``Taylor coefficients''
$m_n:A^{\otimes n}\to A[2-n], n\ge 1$. The condition $d^2=0$ can be rewritten
as a sequence of quadratic equations

$$\sum_{i+j=n+1}\sum_{0\le l\le i}\epsilon(l,j)
m_i(a_0,...,a_{l-1},m_j(a_l,...,a_{l+j-1}),a_{l+j},...,a_n)=0$$
where $a_m\in A$, and $\epsilon(l,j)=
(-1)^{j\sum_{0\le s\le l-1}deg(a_s)+l(j-1)+j(i-1)}$.
In particular, $m_1^2=0$.

\begin{dfn} A morphism of non-unital $\A$-algebras ($\A$-morphism for short)
$(V,d_V)\to (W,d_W)$ is a morphism of tensor coalgebras
$T_+(V[1])\to T_+(W[1])$ of degree zero, which commutes with 
the codifferentials.

\end{dfn}

A morphism $f$ of non-unital $\A$-algebras is determined by its ``Taylor
coefficients'' $f_n:V^{\otimes n}\to W[1-n], n\ge 1$ satisfying
the system of equations

$ \sum_{1\le l_1<...,<l_i=n}(-1)^{\gamma_i}m_i^W(f_{l_1}(a_1,...,a_{l_1}),$\\
$f_{l_2-l_1}(a_{l_1+1},...,a_{l_2}),...,f_{n-l_{i-1}}(a_{n-l_{i-1}+1},...,a_n))=$

$\sum_{s+r=n+1}\sum_{1\le j\le s}(-1)^{\epsilon_s} f_s(a_1,...,a_{j-1},m_r^V(a_j,...,a_{j+r-1}),
a_{j+r},...,a_n).$

Here 
$\epsilon_s=r\sum_{1\le p\le j-1}deg(a_p)+j-1+r(s-j)$,
$\gamma_i=\sum_{1\le p\le i-1}(i-p)(l_p-l_{p-1}-1)+
\sum_{1\le p\le i-1}\nu(l_p)\sum_{l_{p-1}+1\le q\le l_p}deg(a_q)$,
where we use the notation $\nu(l_p)=\sum_{p+1\le m\le i}(1-l_m+l_{m-1})$,
and set $l_0=0$.

\begin{dfn} A non-unital $\A$-category ${\cal C}$
over $k$ is given by the following
data:

1) A class of objects $Ob({\cal C})$.

2) For any two objects $X_1$ and $X_2$ a ${\Z}$-graded $k$-module
of morphisms $Hom(X_1,X_2)$.

3) For any sequence of objects $X_0,...,X_n$, $n\ge 1$,
a morphism of $k$-modules (called a composition map)
$m_n:\otimes_{0\le i\le n-1}Hom(X_i,X_{i+1})\to Hom(X_0,X_n)[2-n]$.

It is required that for any sequence of objects $X_0,...,X_N$, $N\ge 0$
the graded $k$-module $A=A(X_0,...,X_N):=\oplus_{i,j}Hom(X_i,X_j)$,
equipped with the direct sum of the
compositions $m_n, n\ge 1$,
is a non-unital $\A$-algebra.

\end{dfn}

The class of objects $Ob({\cal C})$ will be often denoted
by ${\cal C}$. We hope it will
not lead to a confusion.

\begin{rmk}
A non-unital $\A$-algebra $A$ can be considered as a non-unital
$\A$-category with one object $X$ such that $Hom(X,X)=A$.

\end{rmk}

\begin{dfn} A functor $F:{\cal C}_1\to {\cal C}_2$ between non-unital
$\A$-categories is given by the following data:

1) A map of classes of objects $\phi:{\cal C}_1\to {\cal C}_2$.

2) For any finite sequence of objects $X_0,...,X_n$, $n\ge 0$,
a morphism of graded $k$-modules
$f_n:\otimes_{0\le i\le n-1}Hom_{{\cal C}_1}(X_i,X_{i+1})
\to Hom_{{\cal C}_2}(\phi(X_0),\phi(X_n))[1-n].$

The following condition holds for any $X_1,...,X_N\in {\cal C}_1$:
 the sequence $f_n, n\ge 1$ defines an $\A$-morphism

$$\oplus_{i,j}Hom_{{\cal C}_1}(X_i,X_j)\to \oplus_{i,j}
Hom_{{\cal C}_2}(\phi(X_i),\phi(X_j)).$$

\end{dfn}

\begin{rmk} Let ${\cal C}$ be a non-unital $\A$-category. Let us replace
spaces of morphisms by their cohomology with respect to $m_1$. In other
words, we define $Hom_{H({\cal C})}(X,Y):=\{Ker\, m_1\}/\{Im\, m_1\}$,
where $m_1:Hom_{\cal C}(X,Y)\to Hom_{\cal C}(X,Y)[1]$ is the composition
map. Then $H({\cal C})=({\cal C}, Hom_{H({\cal C})}(\cdot,\cdot))$ gives rise to
a ``non-unital'' category structure with the class
of objects ${\cal C}$ and composition
of morphisms induced by $m_2$. We write ``non-unital'' because there are
no identity morphisms $id_X\in Hom_{H({\cal C})}(X,X)$.

\end{rmk}

\subsection{$\A$-pre-categories}

We start with the notion of {\it non-unital} $\A$-pre-category.
It allows us to work with ``transversal'' sequences of objects.
\footnote{ The notion of ``transversality'' is purely formal in this section.
  The choice of the name will become clear after concrete applications
   in the geometric context, see next sections.}
   Then we will introduce the notion of $\A$-pre-category.
It provides us with a replacement of the identity morphisms.
Roughly speaking, we will have the identity morphism up to
homotopy.

\begin{dfn} Let $k$ be a $\Z$-graded commutative associative ring as before.
A non-unital $\A$-pre-category over $k$ 
 is defined by the following data:

a) A class of objects ${\cal C}$.

b) For any $n\ge 1$ a subclass ${\cal C}_n^{tr}$ of ${\cal C}^n$, 
${\cal C}_1^{tr}={\cal C}$,
called the class of transversal sequences.

c) For $(X_1,X_2)\in {\cal C}_2^{tr}$ a ${\Z}$-graded $k$-module
of morphisms $Hom(X_1,X_2)$.

d) For a transversal sequence of objects $(X_0,...,X_n)$, $n\ge 0$,
a morphism of $k$-modules (composition map)
$m_n:\otimes_{0\le i\le n-1}Hom(X_i,X_{i+1})\to Hom(X_0,X_n)[2-n]$.

It is required that a subsequence $(X_{i_1},...,X_{i_l}),
i_1<i_2<...<i_l$ of a
transversal sequence $(X_1,...,X_n)$ is transversal, and
 that the composition maps satisfy the same system 
of equations as for non-unital $\A$-categories.
Explicitly:

$\sum_{i+j=n+1}\sum_{0\le l\le i}\epsilon(l,j)
m_i(a_0,...,a_{l-1},m_j(a_l,...,a_{l+j}),a_{l+j+1},...,a_n)=0$,\\
where $a_m\in Hom(X_m,X_{m+1})$, and $\epsilon(l,j)=
(-1)^{j\sum_{0\le s\le l-1}deg(a_s)+l(j-1)+j(i-1)}$.

\end{dfn}

\begin{dfn} A functor $F:{\cal C}\to {\cal D}$
 between non-unital
$\A$-pre-categories is given by the following data:

1) A map of classes of objects $\phi:{\cal C}\to {\cal D}$,
such that $\phi^n({\cal C}_n^{tr})\subset {\cal D}_n^{tr}$.

2) For any transversal sequence of objects $(X_0,...,X_n),\,n\ge 1$ in 
${\cal C}$,
a morphism of graded $k$-modules
$$f_n:\otimes_{0\le i\le n-1}Hom_{\cal C}(X_i,X_{i+1})
\to Hom_{\cal D}(\phi(X_0),\phi(X_n))[1-n].$$ 

These data satisfy the following property:
the sequence $f_n, n\ge 1$ defines an $\A$-morphism
$\oplus_{i<j}Hom_{{\cal C}}(X_i,X_j)\to \oplus_{i<j}
Hom_{{\cal D}}(\phi(X_i),\phi(X_j))$.
\end{dfn}

The reader have noticed that we use the summation only 
over the increasing pairs of
 indices $i<j$.
It differs from the case of non-unital $\A$-pre-categories.
The reason is that we do not require the transversality to be
a symmetric relation on objects. It is possible that $Hom(X_0,X_1)$ exists,
but $Hom(X_1,X_0)$ does not. In the case when all $Hom's$ are defined,
two discussed definitions agree.
In particular, a non-unital $\A$-category is the same
as a non-unital $\A$-pre-category such that 
${\cal C}_n^{tr}={\cal C}^n$ for any $n\ge 1$.

\begin{dfn} Let ${\cal C}$ be a non-unital $\A$-pre-category,
$(X_1,X_2)\in {\cal C}_2^{tr}$. We say that $f\in Hom^0(X_1,X_2)$
(zero stands for degree) is a 
quasi-isomorphism if $m_1(f)=0$, and for any objects $X_0$ and $X_3$ such that
$(X_0,X_1,X_2)\in {\cal C}_3^{tr}$ and $(X_1,X_2,X_3)\in {\cal C}_3^{tr}$
one has: $m_2(f,\cdot):Hom(X_0,X_1)\to Hom(X_0,X_2)$ and
$m_2(\cdot, f):Hom(X_2,X_3)\to Hom(X_1,X_3)$ are quasi-isomorphisms
of complexes.

\end{dfn}

\begin{dfn} An $\A$-pre-category is a non-unital $\A$-pre-category ${\cal C}$,
satisfying the following extension property:

For any
finite collection of transversal sequences $(S_i)_{i\in I}$ in ${\cal C}$ and 
an object $X$
there exist objects $X_+$ and $X_-$ and 
quasi-isomorphisms $f_-:X_-\to X$,
$f_+:X\to X_+$  
such that
extended sequences $(X_-,S_i, X_+)$ are transversal for any $
i\in I$.   
\end{dfn}

\begin{rmk} Let ${\cal C}$ be an $\A$-pre-category.
Then partially defined on $H({\cal C})=
({\cal C}, Hom_{H({\cal C})}(\cdot,\cdot)$) composition $m_2$ extends
uniquely, so that it defines a structure of a
category on $H({\cal C})$.

\end{rmk}

\begin{dfn} Let ${\cal C}$ and ${\cal D}$ be $\A$-pre-categories
over $k$. An $\A$-functor $F:{\cal C}\to {\cal D}$ is a functor between
the corresponding non-unital $\A$-pre-categories such that
$F$ takes quasi-isomorphisms in ${\cal C}$ to quasi-isomorphisms
in ${\cal D}$.

\end{dfn}

There is an important notion of equivalence of $\A$-pre-categories
(and $\A$-categories). We are planning to provide
all the details elsewhere (see [KoS]). For the purposes of present paper
we will be using the following definition
(which is in fact a theorem in the more general framework).

\begin{dfn} An $\A$-functor $F:{\cal C}\to {\cal D}$ between
$\A$-pre-categories is called an $\A$-equivalence functor if:

a) Every object $Y\in {\cal D}$ is quasi-isomorphic to an object $\phi(X),
X\in {\cal C}$.

b) The functor induces quasi-isomorphisms of non-unital $\A$-algebras of
morphisms, corresponding to all transversal sequences of objects.

\end{dfn} 

\begin{dfn} Two $\A$-pre-categories ${\cal C}$ and ${\cal D}$ are called 
equivalent if there exists a finite sequence of $\A$-pre-categories 
$({\cal C}_0,\dots,{\cal C}_n),\,{\cal C}_0={\cal C},\,
{\cal C}_0={\cal D}$ such that for every $i,\,0\le i\le k-1$
there exists an $\A$-equivalence functor from ${\cal C}_i$ to 
${\cal C}_{i+1}$ or vice versa.
\end{dfn}

We suggest the language of $\A$-pre-categories in order to replace
more conventional $\A$-categories with strict identity morphisms.

\begin{dfn} An $\A$-category with strict identity morphisms is a
non-unital $\A$-category ${\cal C}$,
such that for any object $X$ there exists an element
$1=1_X\in Hom^0(X,X)$ (identity morphism) such that $m_2(1,f)=m_2(f,1)=f$ and
$m_n(f_1,...,1,...,f_n)=0, n\ne 2$ for any morphisms $f,f_1,...,f_n$.

\end{dfn}

An $\A$-category ${\cal C}$
 with strict identity morphisms is an $\A$-pre-category,
because (in the previous notation)
we can extend a transversal sequence $S$ to
$(X,S,X)$, and set $X_+=X_-=X$, $f_{\pm}=1_X$. Another remark is that
if ${\cal C}$ has only one object, it is an
$\A$-algebra with the strict unit. One can try to
develop the deformation theory of such algebras along the
lines of [KoS1]. The problem is that the corresponding operad
is not free, and the standard theory becomes complicated.
We hope that the framework of $\A$-pre-categories is appropriate
for the purposes of deformation theory of $\A$-categories.
The following conjecture gives another evidence in favor
of such a generalization of $\A$-categories.

\begin{conj} Let us define the notion of equivalent
 $\A$-categories with strict identity
morphisms similarly to the case of $\A$-pre-categories (see above). 
Then the equivalence classes of $\A$-pre-categories are in
one-to-one correspondence with the equivalence classes
of $\A$-categories with strict identity morphisms.

\end{conj}

\subsection {Example: directed $\A$-pre-categories}

There is a useful special case of the notion of $\A$-pre-category
(independently a similar notion was suggested in [Se]).

\begin{dfn} A directed $\A$-pre-category is an  $\A$-pre-category
such that

a) There is bijection of the class of objects
and the set integer numbers: ${\cal C}\simeq {\Z}$. We denote by $X_i$ 
the object
corresponding to $i\in \Z$.

b) Transversal sequences are $(X_{i_1},..., X_{i_n}), i_1<i_2<...<i_n$.

\end{dfn}

The extension property is equivalent to the following one:
 for any object $X_i$
there are exist objects $X_j, j<i$ and $X_m, m>i$ which are quasi-isomorphic
to $X_i$.
Then one can formulate the following version of the previous
conjecture.

\begin{conj} Equivalence classes of directed $\A$-pre-categories
are in one-to-one correspondence with the equivalences classes
of $\A$-categories with strict identity morphisms
and countable
class of objects.

\end{conj}

 Having an $\A$-category ${\cal C}$ with strict identity morphisms,
 and 
countable class of objects, one can construct an infinite
sequence of objects $(X_i)_{i\in {\Z}}$ such that each objects appears
infinitely many times for positive and negative $i$. Then a directed
$\A$-pre-category ${\cal C}^{\prime}$ is defined by setting
$Hom_{{\cal C}^{\prime}}(X_i,X_j)=Hom_{{\cal C}}(X_i,X_j)$ for
$i<j$. All other $Hom's$ are not defined.

\section{ Fukaya category and its degeneration}

\subsection{Fukaya category}

Fukaya category (of a compact symplectic manifold) in the approach presented
here will be in
fact an $\A$-pre-category. 
Our definition is not given in the maximal generality, but it
 will be sufficient for the main application
to abelian varieties. For more elaborated definitions
see [Fu1], [FuOOO].

Let $(V,\omega)$ be a compact symplectic manifold
of dimension $2n$, such that $c_1(T_V)=0\in H^2(V,{\Z})$. 
The Fukaya category
(with the trivial $B$-field) associated with $(V,\omega)$ 
 depends on some additional data, which
we are going to describe below.

We fix an
almost complex structure $J$ compatible with $\omega$ 
and a smooth everywhere non-vanishing
differential form $\Omega$, which is $(n,0)$-form with
respect to $J$. 
Let $L$ be an oriented Lagrangian submanifold. 
Then one
has a map $Arg_L:=Arg_{\Omega_{|L}}:L\to {\R}/2\pi {\Z}$, 
where $Arg_{\Omega_{|L}}(x)$
is the argument
of the non-zero complex number $\Omega(e_1\wedge...\wedge e_n)$,
and $e_1,...,e_n$ is an oriented basis of $T_xL, x\in L$.

\begin{dfn}
Objects of the Fukaya category $F(V,\omega,J,\Omega)$ 
are triples\\
$(L,\rho,\widetilde{Arg}_L)$, where
$L$ is a compact oriented Lagrangian submanifold of $V$ with a spin structure 
(called the {\it support}
of the object), $\rho$ is a 
local system on $L$ (i.e. a complex vector
bundle with flat connection), 
and $\widetilde{Arg}_L:L\to {\R}$ a continuous
lift of $Arg_L$.

We require that for any element
$\beta\in \pi_2^{free}(V,L):=\pi_0(Maps((D^2,\partial D^2), (V,L))$, the 
pairing $([\omega],\beta)$ is equal to zero.

\end{dfn}
 
We will sometimes denote the Fukaya category by $F(V,\omega)$,
or simply by $F(V)$. We will also often
omit from the notation the lifted argument function,
thus denoting an object simply by $(L,\rho)$. The choice of a spin structure
 is essential for  signs, and the choice of the lift 
$\widetilde{Arg}_L$ is necessary for $\Z$-grading.

Let
${\C}_{\varepsilon}$ be the field consisting of formal series
$f=\sum_{i\ge 0}c_ie^{-\lambda_i/\varepsilon}$, 
such that $c_i\in {\C},\lambda_i\in {\R}, \lambda_0<\lambda_1<...,
\lambda_i \to +\infty$. In the case when $[\omega]\in H^2(V,{\Z})$,
one can in fact work over the field ${\C}((q))$, where
$q=exp(-{1\over {\varepsilon}})$.  
In general we equip ${\C}_{\varepsilon}$
with the adic topology: a fundamental system of neighborhoods
of zero consists of sets $U_x=
\{f=\sum_{i\ge 0}c_ie^{-\lambda_i/\varepsilon}|\lambda_i\ge x, i\ge 0\}, x\in {\R} $.
 
 \begin{dfn} For two objects with transversal supports
 we define the space of morphisms such as follows

$$Hom_{F(V,\omega)}((L_1,\rho_1,\widetilde{Arg}_1),
(L_2,\rho_2,\widetilde{Arg}_2)):=(\oplus_{x\in L_1\cap L_2}
Hom(\rho_{1x},\rho_{2x}))\otimes {\C}_{\varepsilon}.$$ 

\end{dfn}

Thus morphisms form
a finite-dimensional 
vector space over the field ${\C}_{\varepsilon}$.
There is a $\Z$-grading of the space of morphisms given
in terms of Maslov index $deg: L_1\cap L_2 \to {\Z}$ (see [Fu2], [Ko], [Se]).
Maslov index depends on a choice of the lift $\widetilde{Arg}_L$
(it corresponds to a choice of a point in the universal covering
of the bundle of the Lagrangian Grassmannians).

\begin{rmk} The condition $([\omega],\beta)=0$ is introduced
for convenience only. It helps to avoid the problem
with the composition $m_0$ we mentioned before. The condition
holds in the case when $V$ is a torus
with the constant symplectic form,
and $L$ is a Lagrangian subtorus. This is our main application
in present paper. In general there is a way to work with
non-trivial $m_0$, if it is small in the adic topology.

\end{rmk}

Now we are going to describe
the $\A$-structure. It is defined by means of a collection of
maps (higher compositions) of graded vector spaces
$m_k^{F(V)}:\otimes_{0\le i\le k-1}
Hom_{F(V)}((L_i,\rho_i),(L_{i+1},\rho_{i+1}))\to 
Hom_{F(V)}((L_0,\rho_0),(L_k,\rho_k))[2-k]$,
where $k\ge 1$ and 
the sequence $(L_0,...,L_k)$ corresponds to a
transversal sequence of objects (the latter notion
will be defined below). 

In the case, when all local systems are trivial of rank one,
the map $m_k$ is defined such as follows. 
Let $D$ be a standard disc $D=\{z\in {\C}| \, |z|\le 1\}$. 
Let us fix a sequence
$(L_0,...,L_k)$ of supports of objects
with pairwise transversal intersections, 
intersection points $x_i\in L_i\cap L_{i+1}, 0\le i\le k-1$,
$x_k\in L_0\cap L_k$, and $\beta \in \pi_2^{free}(V,\cup_{0\le i\le k}L_i)$.
We denote by ${\cal M}(L_0,...,L_k;x_0,...,x_k;\beta)$ the set
of collections $(y_0,...,y_k;\psi)$, where $y_i, 0\le i\le k$ are
cyclically ordered
pairwise distinct points on the boundary $\partial D$, and
$\psi:D\to (V,J)$ a pseudo-holomorphic map such that $\psi(y_i)=x_i,
\psi(\overline{y_iy_{i+1}})\subset L_i, 0\le i\le k, y_0=y_k$,
$[\phi]=\beta$. Here $\overline{y_iy_{i+1}}$ denotes the arc
between $y_i$ and $y_{i+1}$. 
There is a natural action of $PSL(2,{\R})$ on 
${\cal M}(L_0,...,L_k;x_0,...,x_k;\beta)$ arising from the holomorphic
action on $D$ by fractional linear transformations. The action is free
except of the case $k=1, x_0=x_1, \beta=0$, which is not relevant
for our purposes. 

Let $x_i\in L_i\cap L_{i+1}, 0\le i\le k-1, x_k\in L_0\cap L_k$ satisfy
the condition $deg\,x_k=
\sum_{0\le i\le k-1}deg\,x_i+2-k$.
Then the matrix element $(m_k(x_0,x_1,...,x_{k-1}),x_k)$ is given by
the formula
$(m_k(x_0,x_1,...,x_{k-1}),x_k)=\sum \pm q^{({\beta},[\omega])}$,
where sum is taken over 
all $PSL(2,{\R})$-orbits of points in ${\cal M}(L_0,...,L_k;x_0,...,x_k;\beta)$.
Signs are derived
from orientations of certain cycles in the moduli
space ${\cal M}={\cal M}(L_0,...,L_k;x_0,...,x_k;\beta)/PSL(2,{\R})$.
We will comment on them below (see [Fu1], [FuOOO] for more details).
 In the case of non-trivial local systems there is an additional
factor for each summand. It corresponds to the holonomies of 
local system along the arcs.

Now we will describe the transversality condition. Assume that we
are given a sequence of objects $(L_i,\rho_i), 0\le i\le k$ of the Fukaya
category. We say that they are transversal if the following conditions hold:

1) There are only pairwise intersections $L_i\cap L_j$, and they
are transversal.

2) For any subsequence $(L_{i_0},...,L_{i_m}), \,m\ge 1,\,i_0<i_1<...<i_m$, 
any choice of intersection points
 $x_{i_m}\in L_{i_0}\cap L_{i_{m}}$,
$x_{i_p}\in L_{i_p}\cap L_{i_{p+1}}, 0\le p\le m-1$ such that
$deg\,x_{i_m}-
(\sum_{0\le p\le m-1}deg\,x_{i_p}+2-m)=0$, and any $\beta\in
 \pi_2^{free}(V,\cup_{0\le p\le m}L_{i_p})$,
the corresponding component of the moduli space\\
${\cal M}(L_{i_0},...,L_{i_m};x_{i_0},...,x_{i_m};\beta)/PSL(2,{\R})$
contains only smooth points, and is zero-dimensional. 

3) If $deg\,x_{i_m}-
(\sum_{0\le p\le m-1}deg\,x_{i_p}+2-m)<0$ then the corresponding component
is empty.

Let us comment on these conditions (for more details see [FuOOO]).
The first one is needed
to define morphisms. 
The quotient set 
$${\cal M}={\cal M}(L_{i_0},...,L_{i_m};x_{i_0},...,x_{i_m};\beta)/PSL(2,{\R})$$
which appears in the second condition
locally can be identified
with the space of solutions of a non-linear elliptic 
problem. For the linearized
problem the corresponding Fredholm operator has index 
$deg\,x_{i_m}-(\sum_{0\le i\le m-1}deg\,x_{i_p}+2-m)$.
We define smooth points ${\cal M}^{sm}$ of 
${\cal M}$
as such points where the cokernel of the Fredholm operator is trivial.
Then ${\cal M}^{sm}$ is a smooth manifold of the dimension equal to the index.
Moreover, one checks that 
the spaces ${\cal M}^{sm}$ carry natural orientations
given by the determinants of the corresponding Fredholm operators.
 It is here where the choice of spin structures on $L_i$ enters into the game. 
It follows that in the zero-dimensional case what we get is
a set of points with multiplicities $\pm 1$
(in particular, the multiplicities are integer numbers). Multiple
covers and stable maps which appear in the definition of Gromov-Witten
invariants and produce non-trivial denominators, do not appear
in our framework for the Fukaya category (in our case all ``stable maps"
in components of the virtual degree zero 
are embeddings. Ramified coverings appear in higher codimensions).
Therefore one can
define the Fukaya category over 
the ring ${\Z}_{\varepsilon}$ (the integral
version of ${\C}_{\varepsilon}$). The number of points counted with signs
gives a tensor coefficient of $m_k$.

Composition maps satisfy a system
of quadratic equations, thus making $F(V,\omega)$ into a
non-unital $\A$-pre-category.
One can check that it is in fact an $\A$-pre-category.
Proof of the extension property
is based on the following result of Fukaya (see [Fu2], [Se]).

\begin{prp} Let $(L_t,\rho_t)$ be an object obtained
by a small Hamiltonian deformation
of an object $(L,\rho)$ of $F(V)$. Then $(L_t,\rho_t)$
and $(L,\rho)$ are quasi-isomorphic.

\end{prp}

For example, a sequence consisting of one object $(L,\rho)$
can be extended to a transversal
sequence $((L_{t_1},\rho_{t_1}),(L,\rho))$. Similarly, one can extend
any finite set of transversal sequences.

It is easy to see that the set of connected
components of the space of pairs $(J,\Omega)$ (equipped with
the natural topology) is a principal homogeneous space
over the lattice $H^1(V,{\Z})$. Namely, $f:V\to U(1)$ acts
on $(J,\Omega)$ such as follows: $(J,\Omega)\mapsto (J,f\Omega)$.
The following theorem can be derived from [Fu2].

\begin{thm} There exists a set $\Sigma$ of the second category
(in the sense of Baire) in the space of almost complex structures
compatible with $\omega$ such that Fukaya categories
$F(V,\omega,J_1,\Omega_1)$ and $F(V,\omega,J_2,\Omega_2)$
are equivalent as long as $J_1,J_2\in \Sigma$, and
$(J_1,\Omega_1)$ is homotopic to $(J_2,\Omega_2)$.

\end{thm}

Therefore the equivalence class of the Fukaya category depends
on the connected component of the space of pairs.

\begin{rmk} Definitions of this section (and other sections of the paper)
can be modified in order to accomodate the case of non-zero
$B$-field. We are not going to do that in order to avoid
more complicate notations.

\end{rmk}

\subsection{Fukaya-Oh category for torus fibration}

Let $(Y,g_Y,\nabla)$ be an AK-manifold with integral affine structure.
The covariant lattice is denoted by $T_Y^{\Z}$, as before.
 From now on we will assume that $Y$ is {\it compact}.
This is a severe restriction. It was proven in [CY] that in this case a finite cover of space $Y$
 is a torus with the standard affine structure. It appears
in the collapse of complex abelian varieties.

The manifold $X^{\vee}=T_Y^{\ast}/(T_Y^{\Z})^{\vee}$ is the total space of
the torus bundle $p^{\vee}:X^{\vee}\to Y$.
It carries a 
natural symplectic form $\omega=\omega_{X^{\vee}}$ induced from the
standard one on $T^{\ast}Y$.
We endow $X^{\vee}$ with a 1-parameter family of complex structures 
$J_{\eta}, \eta\to 0$ compatible with $\omega$. Indeed, 
the manifold $X_{\eta}^{\vee}:=T_Y^{\ast}/\eta (T_Y^{\Z})^{\vee}$ carries 
a canonical complex structure described before.
We identify  $X^{\vee}$ and $X_{\eta}^{\vee}$ by the map
$(y,v)\mapsto (y,\eta v)$, where $y\in Y, v\in T^{\ast}_{Y,y}$.
Using this identification, we pull back to $X^{\vee}$ the 
complex structure and the metric.
The fibers of $p^{\vee}:X^{\vee}\to Y$ are flat Lagrangian tori for all 
values of $\eta$.

We define on $(X^{\vee}, J_{\eta})$ a
nowhere vanishing $(n,0)$-form $\Omega_{\eta}$ such as follows.
Let us fix an oriented orthonormal basis $e_1,...,e_n$ in
$T_{Y,y}^{\ast}, y\in Y$. We define $\Omega_{\eta}$
as the $n$-form on $X^{\vee}$, which is invariant with respect to the 
$T_{Y,y}^{\ast}/(T_{Y,y}^{\Z})^{\vee}$-action, and is
equal to $\bigwedge_{1\le j\le n}((p^{\vee})^{\ast}e_j+
\sqrt{-1}J_{\eta}(p^{\vee})^{\ast}e_j)$.

Let $L$ be a compact oriented Lagrangian submanifold of $X^{\vee}$ such that
$p^{\vee}_{|L}$ is an unramified covering, and the orientation
of $L$ is induced from the orientation of $Y$. We claim that there is
a canonical choice $\widetilde{Arg}_L^{can}:L\to {\R}$
for the function $\widetilde{Arg}_L:L\to {\R}$.
Indeed, for any point $x\in X^{\vee}$ the space of    Lagrangian
subspaces in $T_{X^{\vee},x}$, which are transversal to the
vertical tangent space $T_x^{vert}=Ker(p^{\vee})_{\ast}$ is contractible.
Let us consider the space ${\cal L}$ of pairs $(x,l)$ such that
$x\in X^{\vee}$ and $l\subset T_{X^{\vee},x}$
is a Lagrangian subspace, which is transversal to $T_x^{vert}$,
and endowed with the orientation induced from $Y$.
Then the function 
$(x,l)\mapsto Arg_{(\Omega_{\eta})_{|l}}(x)\in {\R}/2\pi {\Z}$
admits a unique continuous lifting $\widetilde{Arg}:{\cal L}\to {\R}$,
vanishing at $(x,J_{\eta}(T_x^{vert})),
x\in X^{\vee}$.
Restricting this function to $L$ we obtain $\widetilde{Arg}_L^{can}$.

We will denote by $F^{\eta}(X^{\vee})$ 
the Fukaya category  
 $F(X^{\vee},\omega,J_{\eta},\Omega_{\eta})$, and
by $F_{unram}^{\eta}(X^{\vee})$ its full
$\A$-pre-subcategory with objects  $(L,\rho,\widetilde{Arg}_L^{can})$ such that
$L$ is a compact Lagrangian submanifold with the orientation induced from $Y$,
$p^{\vee}_{|L}$ is an unramified
covering, and $\widetilde{Arg}_L^{can}$ 
was described above. To simplify the notations we will denote
objects of these categories by $(L,\rho)$.

\begin{rmk}
One can check that for transversal 
Lagrangian submanifolds $L_1$ and $L_2$ as above, the Maslov index
at any $x\in L_1\cap L_2$ is equal to the Morse index
at $p^{\vee}(x)$ of the smooth Morse function $f_1-f_2:Y\to {\R}$
such that locally near $x$ one has $L_i=graph\,(df_i)\, (mod(T_Y^{\Z})^{\vee}), i=1,2$.

\end{rmk}

It follows from the results of [FuO] that there exists a limit
of the family of $\A$-pre-categories
$F_{unram}^{\eta}(X^{\vee})$, $\eta\to 0$ in the following sense.
Objects and morphisms of $F_{unram}^{\eta}(X^{\vee})$
do not depend on $\eta$ and remain the same in the limit.
The compositions $m_k^{F_{unram}^{\eta}(X^{\vee})}$ have limits as $\eta\to 0$
in the adic topology of ${\C}_{\varepsilon}$. They will be 
explicitly described below.

The following result can be derived from [FuO].

\begin{prp}
The limiting $\A$-pre-category 
 is equivalent to $F_{unram}^{\eta}(X^{\vee})$
 for all sufficiently small $\eta$.
 
\end{prp}
 
We will denote this $A_{\infty}$-pre-category by $FO(X^{\vee})$
and call it the {\it Fukaya-Oh category} of $X^{\vee}$
(or {\it degenerate Fukaya category} of $X^{\vee}$).

\begin{rmk} In what follows we will assume that $dim\,Y>1$.
The case $dim\,Y=1$ is somewhat different, but also it is much more simple
(see for example [P1]). In particular, $F_{unram}^{\eta}(X^{\vee})$
does not depend on $\eta$ in this case.

\end{rmk}

As we said before,
the objects and morphisms for $FO(X^{\vee})$ are the same as for
$F_{unram}^{\eta}(X^{\vee})$. In order to define
the composition map 
$$m_k:\otimes_{0\le i\le k-1}Hom((L_i,\rho_i),
(L_{i+1},\rho_{i+1}))\to Hom((L_0,\rho_0),(L_{k},\rho_{k}))[2-k]$$
one uses the standard formulas, but the sum runs  
 over certain two-dimensional surfaces in $X^{\vee}$ described below. 
 For a sequence $((L_0,\rho_0),...,(L_k,\rho_k)), k\ge 1$ 
of objects in $FO(X^{\vee})$ we consider immersed
two-dimensional surfaces $S\to X^{\vee}$
such that:

a) Boundary of $S$ belongs
 to $L_0\cup...\cup L_k$.

b) $S=(\cup_{\alpha}T_{\alpha})\cup 
(\cup_{\beta}S_{\beta})$ where
and $T_{\alpha}$ are geodesic triangles in fibers of $p^{\vee}$,
hence they are projected to points in $Y$.

c) Each $S_{\beta}$ is a union of
1-parameter families of geodesic intervals contained in fibers of $p^{\vee}$
(i.e. a ``strip'').
Moreover,
$p^{\vee}_{|S_{\beta}}:S_{\beta}\to Y$ is a fibration over a
connected interval   $I_{\beta}:=p^{\vee}(S_{\beta})$ immersed in $Y$.
Fibers of $S_{\beta}$ over the interior points of $I_{\beta}$ are  geodesic intervals
of strictly positive length. Fibers  of $S_{\beta}$ over the boundary points of $I_{\beta}$
are either  edges of triangles $T_{\alpha}$ or intersection points $x_i\in L_i\cap L_{i+1}, 0\le i\le k-1,\,
x_k\in L_0\cap L_k$.

d) Intervals  $I_{\beta}$ are edges of an immersed planar trivalent tree     $\Gamma\subset Y$.
 Points   $p^{\vee}(T_{\alpha})$     are internal vertices of $\Gamma$.     Tail vertices of    $\Gamma$
  are  projections of the intersection points $x_0,\dots,x_k$.

e) Let $r:T^{\ast}_Y\to X^{\vee}$ be the natural fiberwise universal covering.
If the Lagrangian manifolds $r^{-1}(L_i),i=0,...,k$ 
are locally given by 
 differentials of smooth functions $f_i,i=0,...,k$ on $Y$,
 then the edges of $\Gamma$
must be gradient lines of  $f_i-f_j$. Intersection points
of $r^{-1}(L_i)$ and $r^{-1}(L_j)$ correspond to critical points of $f_i-f_j$.

We depict a typical surface below:

\vspace{2mm}

\centerline{\epsfbox{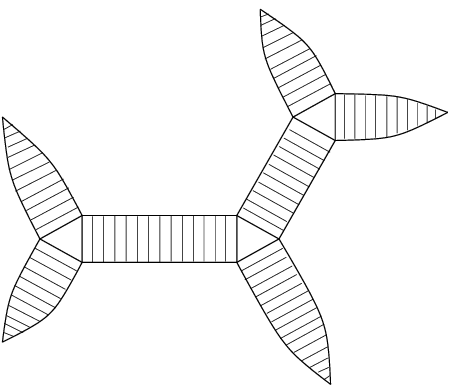}}

The projection of surface $S$ to $Y$ is a gradient tree, with tail vertices
being critical points of $f_i-f_{i+1}$ or of $f_0-f_k$,
and edges $p^{\vee}(S_{\beta})$ being the gradient lines of functions
$f_i-f_j$, where $i=i(\beta), j=j(\beta), i<j$.
The triangles are mapped into the internal vertices of the tree. Here is the picture of
$\Gamma=p^{\vee}(S)$   for surface $S$ as above:

\vspace{2mm}

\centerline{\epsfbox{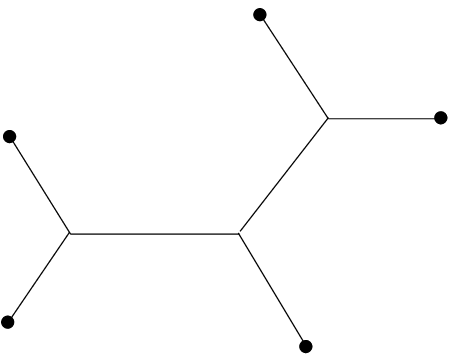}}

Compositions $m_k=m_k^{FO(X^{\vee},\omega)}$ are given
by the standard formulas, but now we are counting
surfaces $S$ described in a)-d). The weight $q^{\langle [S], [\omega]\rangle}$
can be written as $exp(-{1\over{\varepsilon}} \sum_{\beta}var_{p^{\vee}(S_{\beta})}f_{\beta})$,
where $f_{\beta}=f_{i(\beta)}-f_{j(\beta)}$, and {\it var} is
the (positive) variation of the function along the gradient line.

The transversality condition for a sequence of objects of 
Fukaya-Oh category can be formulated similarly to the case
of Fukaya category.

The reader can compare our considerations with those from [FuO].
The fibers of $p^{\vee}$ are ``small'' tori (of the size $O(\eta))$.
The base $Y$ is ``large''  (of the size of $O(1)$).
Hence, the Lagrangian manifolds are close to the zero section of $p^{\vee}$.
This is similar to the situation considered in [FuO].
Indeed, in [FuO] the authors study the $\A$-subcategory of $F(T^{\ast}_Y)$
(where $Y$ is an arbitrary smooth compact manifold),
with the objects $(L,\rho)$ such that $L=\eta\,graph(df)$, $f:Y\to {\R}$
is a smooth function. In other words, they considered Lagrangian sections
of the natural projection $T^{\ast}_Y \to Y$, which are close to the zero
section. 
 When $\eta\to 0$,  
 pseudo-holomorphic discs get ``stretched''
along the fibers of $p^{\vee}$. Thus they look like
the surfaces $S$ described above. Then the higher compositions
of the Fukaya category 
 ``approach'' the compositions
$m_k^{FO(X^{\vee},\omega)}$.
This was proved in [FuO] in the case when $X^{\vee}$ was replaced by
$T^{\ast}_Y$. Considerations from [FuO] apply in our case as well.

\begin{rmk}         One can  extend the Fukaya-Oh category 
considering Lagrangian submanifolds in $X^{\vee}$
 which are not necessarily unramified coverings of $Y$. 
 For example, one can try to add to $FO(X^{\vee})$ new
  objects
  which are local systems on Lagrangian tori which are fibers of the projection  $p^{\vee}:X^{\vee}\to Y$.
 It seems that with these objects one can go much further than with transversal ones. For example,
   in the general case of torus fibrations with singular fibers,  one can argue that for almost any $y\in Y$ there
    are no limiting holomorphic discs with the boundary 
    in the torus $(p^{\vee})^{-1}(y)$ .
    The set of such points $y$ is the complement to a countable union $Z$
     of hypersurfaces in $Y$  (this follows from the fact that
    $dim(Y^{sing})=dim(Y)-2$).   Thus,
      we get a large collection of honest objects without the parasitic 
      composition $m_0$.
      The total picture seems to be quite intricate, as examples show that the subset $Z$ is everywhere dense.
      Presumably, it is related with some mysterious non-abelian
      $1$-cocycle which we will discuss later in the remark  in section 7.1.

\end{rmk}

\section{Morse-Smale complex and the category of Morse functions}

\subsection{Notations from Morse theory}

Let $(Y,g_Y)$ be a compact oriented Riemannian manifold
of dimension $n$, $f:Y\to {\R}$ be
a smooth Morse function.
We will denote the set of critical points of $f$ by
$Cr(f)$. If $x\in Cr(f)$, we will denote by
$U_x$ (resp. $S_x$) the unstable (resp. stable) submanifolds
associated with $x$. Namely, $U_x=\{y\in Y|\,lim_{t\to +\infty}e^{-t\,grad(f)}y=x\}$,
and $S_x=\{y\in Y|\,lim_{t\to +\infty}e^{t\,grad(f)}y=x\}$.

Let $ind(x)$ be the Morse index of $x$, i.e. the negative rank
of the quadratic form $({\partial}^2f)_{|T_xY}, x\in Cr(f)$.
The manifolds $S_x$ and $U_x$ are diffeomorphic to open balls
of dimensions $ind(x)$ and $n-ind(x)$ respectively.
It follows that the
cohomology of $S_x$ with compact support
is a graded vector space with the
only non-zero $1$-dimensional component
in degree $ind(x)$: $H^{\ast}_{c}(S_x)\simeq {\Z}[-ind(x)]$.
A choice of generator defines an orientation of $S_x$.
If the function $f$ satisfies  Morse-Smale transversality condition,
i.e. for any $x,y\in Cr(f)$ the manifolds $S_x$ and $U_y$ intersect
transversally, then $Y=\sqcup_{x\in Cr(f)}S_x$ is a cell decomposition
of $Y$.
The cohomology $H^{\ast}(Y,{\Z})$ can be computed as the cohomology
of the Morse complex $(M^{\ast}(Y,f),\partial)$, with
the components $M^i(Y,f)=\sum_{x\in Cr(f), ind(x)=i}H^{i}_c(S_x)$.
Let us choose orientations of manifolds $S_x$ for all $x\in Cr(f)$.
We endow $U_x$ with the dual orientations. The graded module
$M^i(Y,f)$ can be identified with 
$\oplus_{0\le i\le n}{\Z}^{Cr_i}[-i]$ where $Cr_i$ is the set
of critical points of $f$ of index $i$. The choice of orientations
 gives a basis $([x])_{ x\in Cr(f)}$ of
$M^{\ast}(Y,f)$.

The differential $\partial$ is the standard Morse differential:
$$\partial([x])=\sum_{y\in Cr(f),ind(y)=ind(x)+1}deg(({U}_x\cap
{ S}_y)/{\R})\cdot [y],$$
where $({U}_x\cap {S}_y)/{\R}$ an oriented $0$-dimensional
manifold (a set of points with signs), and $deg(\cdot)\in \Z$ denotes
the total number of points counted with signs. The action of ${\R}$
arises from the natural reparametrization $x\mapsto x+t$ of
the gradient trajectories.

There is also 
a generalization  $M^{\ast}(Y,f,\rho)$ of the Morse complex
for a flat vector bundle $\rho$ (see [BZ], [HL]).

\subsection{Morse $\A$-category of smooth functions}

Here we will  define the Morse category  of smooth functions $M(Y)$
 following [FuO].
 It will be an $\A$-pre-category  over ${\C}$.
Objects of $M(Y)$ are pairs $(f,\rho)$,
where $f:Y\to {\bf R}$ is a smooth function, and $\rho$ is a 
local system of finite-dimensional complex vector spaces
on $Y$. Before defining the transversality
of objects, we will define the transversality of functions.

Suppose we are given a sequence
of smooth functions $(f_0,...,f_k), k\ge 2$ such that all
$f_i-f_j, i\ne j$ are Morse functions, and a sequence of critical points $x_i\in Cr(f_i-f_{i+1}),
0\le i\le k-1, x_k\in Cr(f_0-f_k)$. We will use oriented binary planar
trees in order to describe certain moduli spaces associated with
such  sequences.
 Let us fix a planar trivalent tree $T$ with $k+1$ tails vertices. Among the tail vertices we choose
 one and call it the {\it root} vertex.      Let us orient edges of $T$ along the shortest paths
 towards the root.  Thus, $T$ becomes a binary tree considered as an {\it oriented} tree .
We depict $T$ inside
of the standard unit disc $D\subset {\R}^2$ in such a way that
tail vertices of $T$ belong to $\partial D$, and connected
components of $D\setminus T$ are cyclically numbered from $0$ to $k$
in the clockwise order.    We assume that numbers attached to
 two regions near the root vertex are $0$ and $k$.

 \vspace{2mm}
\centerline{\epsfbox{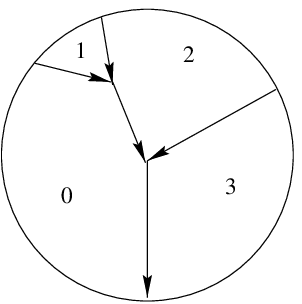}}

We define a gradient immersion of $T$ into $Y$ as a continuous map $j:T\to Y$
such that:

1) The restriction $j|_{e}$ is an orientation preserving homeomorphism 
of the edge $e$
onto an interval in the
gradient line of $f_{l(e)}-f_{r(e)}$,
where the label $l(e)$ (resp. $r(e)$) corresponds
to the region of $D\setminus T$ which is left (resp. right) to $e$.

2) Each tail vertex $v$ is mapped  to the point $x_v\in Cr(f_{l_v}-f_{r_v})$,
where $l_v$ (resp. $r_v$) is the label of the region which is
left (resp. right) to the only tail edge containing $v$.

We will need immersed binary trees (let us call them {\it gradient} trees)
in order to define compositions and transversal sequences in the Morse 
$\A$-pre-category.
These structures can be defined in terms of certain varieties, which
we are going to describe now.

Suppose that we are given a sequence of functions $(f_0,...,f_k), k\ge 2$ and critical points
 $(x_0,\dots,x_k)$
as above, and a binary planar tree $T$.
 Let us consider the
manifold $Y(T)=Y^{V_i(T)}$, where $V_i(T)$ is the set of internal
vertices of $T$. We are going to define several submanifolds in $Y(T)$.
For each tail vertex $v_m, 0\le m \le k-1$ we define 
$Z_{v_m}=\pi_{\hat{v}_m}^{-1}(U_{x_m})$ and for $m=k$
we define $Z_{v_k}=\pi_{\hat{v}_k}^{-1}(S_{x_k})$.
Here $\hat{v}_l$ denotes the second endpoint of the edge of $T$ containing
$v_l$, and $\pi_v:Y(T)\to Y$ is the canonical projection on the factor
corresponding to $v\in V_i(T)$.

    For pair $(f_i, f_j)$ we define a subset $Z_{i,j}\subset Y\times Y$,
consisting of pairs $(y_1,y_2)$ such that $y_1\ne y_2$ and
$y_2=e^{t\,grad(f_i-f_j)}y_1$ for some $t>0$. Then $Z_{i,j}$
is a non-compact submanifold of $Y\times Y$.

 An edge $e$ of $T$ we call internal if both endpoints of it
 are  internal vertices. The set of internal edges
  we denote by $E_i(T)$.
  For each internal edge $e\in E_i(T)$, which separates two regions 
labeled by $l(e)$ (left) and $r(e)$ (right),
 we define 
 a submanifold $Z_{e}=\pi_{e}^{-1}(Z_{l(e),r(e)})$, 
 where $\pi_e:Y(T)\to Y\times Y$
is the natural projection.

It follows from the definitions that the space of gradient immersions of a given $T$
as above, up to homeomorphisms preserving tails, can be identified
with ${\cal M}(T;f_0,...,f_k;x_0,...,x_k):=(\cap_{0\le m\le k}Z_{v_m})
\cap (\cap_{e\in E_i(T)}Z_e)\subset Y^{V_i(T)}$.

\begin{dfn} We say that a sequence $(f_0,\dots,f_k),\,k\ge 2$ 
is $T$-transversal
 for a given tree $T$,       
 if for any sequence of intersection points   $(x_0,\dots,x_k)$ such that
  $$\sum_{i=0}^{k-1} ind(x_i)-ind(x_k)\le k-2$$
 the collection of submanifolds $\bigl(
(Z_{v_m})_{0\le m\le k}, (Z_e)_{e\in E_i(T)}\bigr)$ is transversal
in $Y(T)$ (i.e. intersection of any subcollection 
is transversal).  For $k=1$, we say that $(f_0,f_1)$ is $T$-transversal
(there is only one tree $T$ in this case) if $f_0-f_1$ is a Morse function,
satisfying the Morse-Smale transversality condition.

\end{dfn}

\begin{rmk} As in the case of Fukaya category we consider here only spaces
${\cal M}(T;f_0,...,f_k;x_0,...,x_k)$
of (virtual) dimensions less or equal than zero.   Our condition in the case of strictly negative
  dimension means that the moduli space is empty.

\end{rmk}

It can be proven (see [Fu1]) that there exists a subset  of second
Baire category in $(C^{\infty}(Y))^{\Z}$ such that for any element 
$(f_i)_{i\in\Z}$ of this set and for any strictly
increasing sequence of integers
$i_0<\dots <i_k$ and for any planar tree $T$ with $k+1$ tails,
the sequence
 $(f_{i_0},\dots,f_{i_k})$ is $T$-transversal.

\begin{dfn} A sequence of objects $(f_0,\rho_0),...,(f_k,\rho_k)$
is called transversal if  for any $m\ge 1$ and any binary tree $T$
with $m+1$ tails, 
an arbitrary subsequence $(f_{i_0},...,f_{i_m}), i_0<...<i_m$ is 
$T$-transversal.

\end{dfn}

For any two transversal objects $W_0=(f_0,\rho_0)$ and $W_1=(f_1,\rho_1)$
we define the space of morphisms 
$Hom_{M(Y)}(W_0,W_1)$ as the Morse complex
$M^{\ast}(Y,f_0-f_1, \rho_0^{\ast}\otimes \rho_1)$.
Now we define the 
$A_{\infty}$-structure on $M(Y)$.

The map $m_1:Hom((f_0,\rho_0),(f_1,\rho_1))
\to Hom((f_0,\rho_0),(f_1,\rho_1))[1]$
 is the standard differential in the
Morse-Smale complex. Higher compositions $m_k$ where $k\ge 2$
 for transversal sequences
of objects are linear maps
 $$m_k:\otimes_{0\le i\le k-1}Hom((f_i,\rho_i),(f_{i+1},\rho_{i+1}))
\to Hom((f_0,\rho_0),(f_{k},\rho_{k}))[2-k]$$
Each  $m_k$ is defined as a sum
$m_k=\sum\pm m_{k,T}$ where $T$ runs through the set of isomorphism
classes of oriented binary planar trees with $(k+1)$ tails. Let us describe the summands $m_{k,T}$.
For simplicity we will give the formulas in the case
when all local systems are trivial of rank one.

Let us fix
critical points $y_{i}\in Cr(f_i-f_{i+1}),0\le i\le k-1,y_k\in Cr(f_0-f_k)$,
such that $\sum_{0\le i\le k-1}ind(y_i)=ind(y_k)+2-k$,
and orientations of manifolds $S_{x_i}, 0\le i\le k$.
It follows from the definition of a transversal sequence that  the moduli space
of gradient trees ${\cal M}(T;f_0,...,f_k;y_0,...,y_k)$
is an oriented compact zero-dimensional manifold.

\begin{dfn} We define compositions $m_k, k\ge 2$ by the formula

$$m_k([y_0],...,[y_{k-1}])=\sum_{[T]}\sum_{y_k\in Cr(f_0-f_k)}
deg( {\cal M}(T;f_0,...,f_k;y_0,...,y_k))\cdot [y_k]$$
where $[T]$ is the equivalence class of $T$ as an abstract oriented planar tree,
and $deg(\cdot)\in \Z$ is the total number of points counted with signs, as
before.

\end{dfn}

For local systems of higher ranks one proceeds as
in the case of Fukaya categories, using flat
connections in order to define an analog of the holonomy 
of local systems.

One can obtain slightly different formulas for $m_k$ in the following way.
For any point $\gamma\in {\cal M}(T;f_0,...,f_k;y_0,...,y_k)$
we define the weight

$$ w_{\gamma}=exp(-{1\over{\varepsilon}}\sum_{e\in E(T)}
var_{\gamma}(f_{l(e)}-f_{r(e)})) \in {\C}_{\varepsilon}.$$
Here
$var_{\gamma}(f_{l(e)}-f_{r(e)})>0$ is a variation
of $f_{l(e)}-f_{r(e)}$
along the gradient line $\gamma(e)$, which is defined such as follows:
$var_{\gamma}(f_{l(e)}-f_{r(e)})=(f_{l(e)}-f_{r(e)})(y_{max})-
(f_{l(e)}-f_{r(e)})(y_{min})$,
where $y_{max}$ and $y_{min}$ are the endpoints  of $\gamma(e)$, such that
$(f_{l(e)}-f_{r(e)})(y_{max})-
(f_{l(e)}-f_{r(e)})(y_{min})>0$.
After extension of scalars to ${\C}_{\varepsilon}$ one can choose
another basis in $Hom_{M(Y)}(W_0,W_1)$, namely
$[y]_{new}=[y]exp({(f_0(y)-f_1(y))\over {\varepsilon}})$ for
$y\in Cr(f_0-f_1)$. Then the formulas for $m_k$ will be modified.
The contribution of each $\gamma$ will be multiplied by $w_{\gamma}$.
The formulas will be similar to those for the Fukaya-Oh category (see Section 5.2).

\subsection{De Rham $\A$-category of smooth functions}

The other $A_{\infty}$-pre-category we are interested in
 will be a differential-graded 
category (dg-category for short). In other words, it is an
$\A$-category with strict identity morphisms and vanishing 
compositions $m_n, n\ge 3$. 
We will call it {\it de Rham category}
of $Y$ and denote by $DR(Y)$. 
Objects of $DR(Y)$ are same as for $M(Y)$. They are pairs $(f,\rho)$, 
 where $f: Y\to {\bf R}$ is a smooth function and
 $\rho$ is a local system on $Y$. 
Morphisms are complexes defined by the formula
$$Hom_{DR(Y)}((f_0,\rho_0),(f_1,\rho_1))=
\Gamma(Y, \wedge^{\ast}T^{\ast}_Y\otimes Hom(\rho_0,\rho_1)).$$
Notice that the space
of morphisms does not depend on $f_0$ and $f_1$. 
The composition
of morphisms is defined in the obvious way: in a local trivialization of
$\rho_0$ and $\rho_1$ it is given by the product
of matrices with the coefficients in $\Omega^{\ast}(Y)$.

Now we can formulate the main result of this section.

\begin{thm} $\A$-pre-categories  $M(Y)$ and $DR(Y)$ are equivalent. 
 
\end{thm}

The proof of the theorem will occupy the rest of the section.
First, we will discuss
a version of formulas from ``homological perturbation theory''
(see [GS], [Me]).
They will give an $\A$-structure on a subcomplex of a dg-algebra.
Then we will discuss an approach to the proof based on the ideas of [HL].
It seems plausible that
an alternative proof (but, presumably, much more difficult)
 can be obtained within the framework of Witten
complex, using methods of [BZ].

\subsection{$\A$-structure on a subcomplex}

In this section we are going to restate
in a convenient form some results from [GS] and [Me].

Let $(A, m_n), n\ge 1$ be
a non-unital $A_{\infty}$-algebra, $\Pi:A\to A$ be an idempotent which
commutes with the differential $d=m_1$. In other words, $\Pi$ is a linear
map of degree zero such that
$d\Pi=\Pi d, \Pi^2=\Pi$. Assume that we are given an homotopy
$H:A\to A[-1]$, $1-\Pi=dH+Hd$. Let us denote the image of $\Pi$ by $B$. Then
we have an embedding $i:B\to A$ and a projection $p:A\to B$, such that
$\Pi=i\circ p$.

Let us introduce a sequence of linear operations 
$m_n^B:B^{\otimes n}\to B[2-n]$ in the following way
\footnote{In a special case similar formulas appeared in [Go]}:

a) $m_1^B:=d^B=p\circ m_1\circ i$;

b) $m_2^B=p\circ m_2\circ (i\otimes i)$;

c) $m_n^B=\sum_T \pm m_{n,T}, n\ge 3$.

Here the summation is taken over all oriented planar trees $T$
with $n+1$ tails vertices (including the root vertex), such that the (oriented) valency
$|v|$ (the number of ingoing edges) of every internal vertex of $T$
is at least $2$.
In order to describe the linear map $m_{n,T}:B^{\otimes n}\to B[2-n]$ we need
 to make some preparations.
  Let us consider another tree $\bar{T}$ which is obtained
from $T$ by the insertion of a new vertex into
every internal edge. As a result, there will be two types of internal
vertices in ${\bar{T}}$: the ``old'' vertices, which coincide with the
internal vertices of $T$, and
the ``new'' ones, which can be thought geometrically as the midpoints of the
internal
edges of $T$. 

To every tail vertex of ${\bar{T}}$ we assign the embedding $i$.
To every ``old'' vertex $v$ we assign $m_k$ with
$k=|v|$. To every ``new''
vertex we assign the homotopy operator $H$. To the root we assign the projector
$p$. Then moving along the tree down to the root one 
 reads off the map $m_{n,T}$ as the composition of maps assigned
to vertices of $\bar{T}$. Here is an example of $T$ and
$\bar{T}$:

\vspace{3mm}

\centerline{\epsfbox{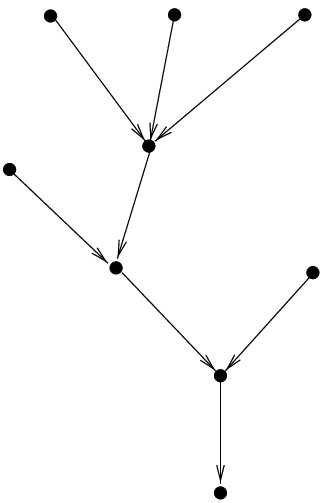}}

\vspace{3mm}

\centerline{\epsfbox{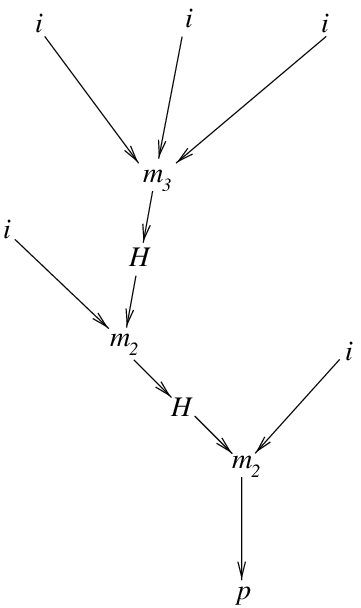}}

\vspace{2mm}

\begin{prp} The linear map $m_1^B$ defines a differential in $B$.

\end{prp}

{\it Proof.} Clear. $\blacksquare$

\begin{thm} The sequence $m_n^B, n\ge 1$ gives rise to a structure of an
$\A$-algebra on $B$.

\end{thm}

{\it Sketch of the proof.} The proof is
quite straightforward, so we just briefly show main 
steps of computations.

First, one observes that $p$ and $i$ are
homomorphisms of complexes. In order to prove the theorem we will
replace for a given $n\ge 2$ each summand $m_{n,T}$ 
by a different one, and then compute
the result in two different ways. Let us consider
a collection of trees $\{\bar{T}_e\}_{e\in E(\bar{T})}$ such that $\bar{T}_e$ 
is obtained from $\bar{T}$ in the following way:

a) we split the edge $e$ into two edges by inserting a new vertex $w_e$
inside $e$;

b) the remaining part of $\bar{T}$ is unchanged.

We assign $d=m_1$ to the vertex $w_e$, and keep 
all other assignments untouched.
In this way we obtain a map $m_{n,\bar{T}_e}:B^{\otimes n}\to B[3-n]$.

Let us consider the following sum (with appropriate signs):

$$\hat{m}_n^B=\sum_T\sum_{e\in E(\bar{T})}\pm m_{n,\bar{T}_e}.$$

We can compute it in two different ways: using the relation $1-\Pi=dH+Hd$,
and using the formulas for $d(m_j),\,j\ge 2$ given by 
the $A_{\infty}$-structure
on $A$. The case of the relation $1-\Pi=dH+Hd=:d(H)$ gives

$$\hat{m}_n^B=d(m_n^B)-m_n^{B,\Pi}+m_n^{B,1}$$
where $m_n^{B,\Pi}$ is defined analogously to $m_n^B$,
with the only difference that we assign to a new vertex operator $\Pi$
instead of $H$ for some edge $e\in E_i(T)$. 
Similarly, the summand $m_n^{B,1}$
is defined if we assign to a new vertex operator $1=id_A$ instead of $H$.
Formulas for $d(m_j)$ are quadratic expressions in $m_l,\,l<j$.
This gives us another identity
$$\hat{m}_n^B=m_n^{B,1}$$
Thus we have $d(m_n^B)=m_n^{B,\Pi}$,
and it is exactly the $\A$-constraint for the collection
$(m_n^B)_{n\ge 1}$. $\blacksquare$

Moreover, using similar technique, one can prove the following result.
                                      
\begin{prp} There is a canonical $\A$-morphism $g:B\to A$, 
which defines a quasi-isomorphism
of $A_{\infty}$-algebras.

\end{prp}

For the convenience fo the reader we give an explicit formula for a
 canonical choice of $g$. The operator $g_1:B\to A$ is defined as the 
 inclusion $i$.
 For $n\ge 2$ we define $g_n$ as the sum of terms $g_{n,T}$ 
 over all planar trees
 $T$ with $n+1$ tails. Each term $g_{n,T}$ is similar to the term
  $m_{n,T}$ defined above,
 the only difference is that 
 we insert operator $H$ instead of $p$ into the root vertex.
 
 One can also construct an explicit $\A$-quasi-isomorphism $A\to B$.

\begin{rmk}

a) Similar construction works
in the case of an arbitrary non-unital $\A$-category. In that case one needs
projectors $\Pi_{X,Y}$ and
homotopies $H_{X,Y}$ 
for every graded space of morphisms $Hom(X,Y)$. All formulas
remain the same as in the case of $\A$-algebras.
The resulting $\A$-category with the spaces of morphisms given
by $\Pi_{X,Y}(Hom(X,Y))$ is equivalent to the original one.
We will use this fact later.

b) Propositions 4 and 5 should hold in a much more general case of algebras
over operads (see e.g. [M]).

\end{rmk}

\subsection{Projectors and homotopies in Morse theory}                                                                  

We would like to apply formulas for the $\A$-structure on a subcomplex
to the proof of the Theorem 2. In order to do that we need to identify
the Morse complex with a direct summand of the de Rham complex.
Our approach is based on the ideas of Harvey and Lawson (see [HL]).

Let $Y$ be a compact oriented smooth manifold, $dim\,Y=n$.
The space of currents $D^{\prime}(Y)$ we will identify with
the space of distribution-valued differential forms.
Continuous linear operators $\Omega^{\ast}(Y)\to D^{\prime}(Y)$
are given by their Schwartz kernels, which are elements
of $D^{\prime}(Y\times Y)$. Smoothening operators
$D^{\prime}(Y)\to \Omega^{\ast}(Y)$ have kernels in 
$\Omega^{\ast}(Y\times Y)\subset D^{\prime}(Y\times Y)$.

With any oriented submanifold $Z\subset Y$ , $\dim\,Z=k$ of finite volume
we associate a canonical current $[Z]$ of degree $n-k$ (namely,
we can integrate smooth $k$-forms over $Z$).

Let $g_Y$ be a Riemannian metric on $Y$, and $f$ be a Morse-Smale function.
The gradient flow $exp(t\,grad(f)), t\ge 0$ gives rise to a
$1$-parameter semigroup acting on $\Omega^{\ast}(Y)$:
$\psi^t(\alpha)=exp(t\,grad(f))_{\ast}(\alpha)$.
Schwartz kernel of $\psi^t$ is $[G_t]$ where manifold $G_t\subset Y\times Y$
is given by $G_t:=graph(exp(t\,grad(f)))$.
We also have the identity

$$id-\psi^t=dH^t+H^td,$$
where $H^t:\Omega^{\ast}(Y)\to \Omega^{\ast}(Y)\subset D^{\prime}(Y)$
is a linear operator of degree $-1$ defined by the distributional
kernel $[Z_t]$, $Z_t:=\cup_{0\le t^{\prime}\le t}graph(exp(t^{\prime}\,grad(f)))$.

It is checked in [HL] that this picture has a limit 
as $t\to +\infty$.
Namely, there exist limits of currents $[G_t]$ and $[Z_t]$:

$$[G_{\infty}]=lim_{t\to +\infty}[G_t]=\sum_{x\in Cr(f)}[S_x]\times [U_x]$$
$$[Z_{\infty}]=lim_{t\to +\infty}[Z_t]=[\cup_{0\le t<+\infty}G_t]$$
Linear operators $\psi^{\infty}$ (of degree zero)
 and $H^{\infty}$ (of degree $-1$), corresponding
to these kernels, map $\Omega^{\ast}(Y)$ to  $D^{\prime}(Y)$
and satisfy the identity

$$i-\psi^{\infty}=dH^{\infty}+H^{\infty}d,$$
where $i:\Omega^{\ast}(Y)\to D^{\prime}(Y)$ is the natural
inclusion. According to the de Rham theorem this inclusion
is a quasi-isomorphism of complexes, therefore $\psi^{\infty}$ is.
Morally, $\Pi_{\infty}:=\psi^{\infty}$ should be thought of as
a projector. The image $\Pi_{\infty}(\Omega^{\ast}(Y))\subset D^{\prime}(Y)$
coincides with $\oplus_{x\in Cr(f)}{\R}\cdot [U_x]$.
We have
$$\Pi_{\infty}(\alpha)=\sum_{x\in Cr(f)}(\int_{S_x}\alpha)\cdot[U_x]=
\sum_{x\in Cr(f)}\int_Y(\alpha \wedge [S_x])\cdot [U_x].$$

Moreover, the operator $\Pi_{\infty}$ commutes with the differentials.
Hence the complex
$\Pi_{\infty}(\Omega^{\ast}(Y))$
is a finite-dimensional subcomplex
of $D^{\prime}(Y)$ isomorphic to the Morse complex $M^{\ast}(Y,f)$. In fact it is
quasi-isomorphic to both complexes $\Omega^{\ast}(Y)$ and $D^{\prime}(Y)$.
In this way Harvey and Lawson prove that the de Rham cohomology is
isomorphic to the cohomology of Morse complex.

In order to construct actual projectors and homotopies we will
proceed as follows.
Let $\rho_{\delta},\delta\to 0$ be a family of smooth closed
differential $n$-forms on $Y\times Y$ such that $supp(\rho_{\delta})$
belongs to the open $\delta$-neighborhood $N_{\delta}$ of the diagonal
$diag\subset Y\times Y$, and the cohomology class of $\rho_{\delta}$
in $H^n_{c}(N_{\delta},{\R})$ is the same as of $[diag]$.

We define $R_{\delta}:D^{\prime}(Y)\to \Omega^{\ast}(Y)$ as the
integral operator given by the kernel $\rho_{\delta}$.

\begin{lmm} 1) The operator $R_{\delta}$ is a homomorphism
of complexes.

2) If $Z_1, Z_2\in Y$ are two oriented submanifolds of finite volume
such that they intersect transversally at finitely many points,
and $dim\,Z_1+dim\,Z_2=dim\,Y$, 
$\overline{Z}_1\cap \overline{Z}_2={Z}_1\cap {Z}_2$, then for
sufficiently small $\delta$ one has:

$$\int_YR_{\delta}([Z_1])\wedge R_{\delta}([Z_2])=deg(Z_1\cap Z_2)\in\Z$$

3) There exists a linear operator $h_{\delta}:\Omega^{\ast}(Y)\to 
\Omega^{\ast}(Y)$ such that its kernel has support in $N_{\delta}$,
the wave front $WF(h_{\delta})$ is the conormal bundle of 
$diag\subset Y\times Y$, and

$$dh_{\delta}+h_{\delta}d=id-(R_{\delta})_{|\Omega^{\ast}(Y)}.$$

\end{lmm}

{\it Proof .} Part 1) follows from the fact that $\rho_{\delta}$
is a closed current. Part 2) follows from the fact that
$R_{\delta}$ changes the supports of $Z_i, i=1,2$ by $O(\delta)$.
To prove part 3) one observes that the operators $id$ and
$(R_{\delta})_{|\Omega^{\ast}(Y)}$ preserve the space of smooth
forms $\Omega^{\ast}(Y)$, and $\rho_{\delta}$ is cohomologous
to $[diag]$. $\blacksquare$

Let $x,y\in Cr(f)$ be two critical points of the same Morse index.
Then $deg(S_x\cap U_y)=\delta_{xy}$ (the Kronecker symbol). By the part 2) of the
Lemma, for sufficiently small $\delta$ we obtain the identity

$$ \int_YR_{\delta}([S_x])\wedge R_{\delta}([U_y])=\delta_{xy}$$
This implies the following result.

\begin{prp} Let us define for a sufficiently small $\delta$
 a linear operator $D^{\prime}(Y)\to \Omega^{\ast}(Y)$
by the formula
$\Pi_{\delta}(\alpha)=\sum_{x\in Cr(f)}(\int_Y\alpha \wedge R_{\delta}([S_x]))
\cdot R_{\delta}([U_x]).$

Then 

1) $\Pi_{\delta}^2(\alpha)=\Pi_{\delta}(\alpha)$ 
if $\alpha\in \Omega^{\ast}(Y)$, and 
$ \Pi_{\delta}d=d\Pi_{\delta}.$

2) The image $\Pi_{\delta}(M^{\ast}(Y,f))$ is a subcomplex
in $\Omega^{\ast}(Y)$ which is canonically isomorphic to the Morse
complex $M^{\ast}(Y,f)$.

\end{prp}

We define a homotopy operator 
$H_{\delta}:\Omega^{\ast}(Y)\to \Omega^{\ast}(Y)[-1]$ as an integral operator
given by the kernel
$(R_{\delta}\boxtimes R_{\delta})[Z_{\infty}]+
(h_{\delta}\boxtimes h_{\delta})([diag]).$
(The last summand is well-defined because of the condition on
the wave front of $h_{\delta}$).
It is easy to check that the following identity holds:

$$id-\Pi_{\delta}=dH_{\delta}+H_{\delta}d .$$

Thus we have a family of homotopies and projectors parametrized by
$\delta$.

\begin{rmk} One can define the projector $\Pi_{\delta}$ using
another canonical element $\sum_{x\in Cr(f)}[S_x]\otimes R_{\delta}([U_x])$,
instead of $\sum_{x\in Cr(f)}R_{\delta}([S_x])\otimes R_{\delta}([U_x])$,
as we did. The above Proposition holds for the new canonical element
as well.

\end{rmk}

There is a version of the previous construction, which will be useful
in the next subsection. Namely, we start with a differential
$n+1$-form $\rho$ on $Y\times Y\times (0,1)$ such that for 
 the support of
 $supp(\rho)$ belongs to $\sqcup_{\delta>0} (N_{\delta},\delta)$ 
for all sufficiently small $\delta\in (0,1)$, and
$\rho$ defines the same cohomology class in $H_c^n(Y\times Y\times (0,1))$
as $[diag]\times (0,1)$. 

Let us consider now the spaces
$\Omega_0^{\ast}(Y):=\varinjlim_{\delta\to 0}\Omega^{\ast}(Y\times (0,\delta))$
and 
$D_0^{\prime}(Y):=\varinjlim_{\delta\to 0}\Omega^{\ast}(0,\delta)\widehat{\otimes} D'(Y).$
It is easy to see that both complexes   $\Omega_0^{\ast}(Y)$ and   $D_0^{\prime}(Y)$
are quasi-isomorphic to  $\Omega^{\ast}(Y)$ .

 We define a linear operator $R:D_0^{\prime}(Y)\to \Omega_0^{\ast}(Y)$
similarly to the definition of $R_{\delta}$. Then the Lemma
and the Proposition hold with obvious changes. We will denote
the corresponding objects by the same letters as before, skipping
the subscript $\delta$ (like $H$ for the homotopy and $\Pi$ for
the projector). 
Morally, they are obtained from the old
objects by extending them as differential forms 
``in the direction of $\delta$''.

\subsection{Proof of the theorem from 6.3}

For simplicity we will assume that all local systems 
are trivial and have rank one. The general case is
completely similar.

We are going to construct the following chain of $\A$-equivalences
 connecting  $DR(Y)$   and
  $M(Y)$:
  $$DR(Y)\mono DR_0(Y)\hookleftarrow DR_0^{tr}(Y)\hookleftarrow
   DR_0^{tr,\Pi}(Y))
 \leftarrow M(Y)$$
  
  Classes of objects of all these categories will be the same, and all functors
  will be identical on objects.
  
  The $\A$-pre-category $DR_0(Y)$ is in fact a dg-category, i.e. all
  sequences of objects are transversal, compositions $m_k$ vanish for $k\ge 3$
  and it has strict identity morphisms.
   The space  
$Hom_{DR_0(Y)}(f_0,f_1)$ is defined as $\varinjlim_{\delta\to 0}\Omega^{\ast}(Y\times (0,\delta))
=\Omega_0^{\ast}(Y).$
Clearly the space of morphisms does not depend on objects.
Using the wedge product of differential forms we make
$DR_0(Y)$ into a dg-category over the field $\C$.
There is a natural functor $DR(Y)\to DR_0(Y)$, which is the
identity map on objects.  On morphisms it is
the natural embedding of $\Omega^{\ast}(Y)$
as the subspace of forms on $Y\times (0,\delta)$,
which are pullbacks of forms on $Y$.
Clearly it establishes an equivalence of $\A$-categories.

The $\A$-pre-category $DR_0^{tr}(Y)$ is defined as the full subcategory
of $DR_0(Y)$, and it differs from the latter only by  the choice
of transversal sequences. Namely, we use the same notion of transversality  in   $DR_0^{tr}(Y)$
 as in the Morse category.
 
The next  $\A$-pre-category $DR_0^{tr,\Pi}(Y)$ is obtained from 
$DR_0^{tr}(Y)$ by applying homological perturbation theory.
 For any two transversal objects $f_0,f_1$ of $DR_0^{tr}(Y)$
we define $Hom_{DR_0^{tr,\Pi}(Y)}(f_0,f_1)$ as $\Pi_{f_0,f_1}(\Omega^{\ast}_0(Y))$.
Here $\Pi_{f_0,f_1}$ is the projector $\Pi$ corresponding to the Morse function $f_0-f_1$,
 it was described at the end of the
previous subsection. 
  We also have
homotopies $H_{f_0,f_1}$ associated with $f_0-f_1$.
Then formulas
of homological perturbation theory (summation over trees)
give rise to an $\A$-pre-category $DR_0^{tr,\Pi}(Y)$ and an equivalence
   $DR_0^{tr,\Pi}(Y)\to DR_0^{tr}(Y)$.

The last functor $\Psi:M(Y)\to DR_0^{tr,\Pi}(Y))$ will 
have no non-trivial higher components $\Psi_n$ 
for $n\ge 2$. The first component $\Psi_1$ of it is a linear map
 $$\Psi_1:Hom_{M(Y)}(f_0,f_1)\to Hom_{DR_0^{tr,\Pi}(Y)}(f_0,f_1)$$
for every transversal pair $(f_0,f_1)$.
 Recall that $Hom_{M(Y)}(f_0,f_1)$ has a basis $\{[x]\}$ labeled by
critical points $x\in Cr(f_0-f_1)$.  We define 
$\Psi_1([x])$ as $R([S_x])$.
It is clear that $\Psi_1$ gives a quasi-isomorphism of complexes
 for every transversal pair $(f_0,f_1)$.

 Now, we claim that $\Psi$ is an $\A$-functor. This means that
 $\Psi_1$ maps all higher compositions in $M(Y)$ to higher
 compositions in $DR_0^{tr,\Pi}(Y)$.
 This follows directly from the descriptions of higher compositions
 in both
  categories in terms of planar trees, see Sections 6.2, 6.4.
Notice that the number of functions
in a given sequence of objects is finite. For all sufficiently small $\delta$
every summand in the formula for $m_k^{M(Y)}$, corresponding
to a binary tree $T$, coincides with the summand for
$m_k^{M(Y)}$ corresponding to the same $T$ 
(we can assume that
 $\delta$ is so small that the part 2) of the Lemma 2 can be applied).
  The theorem is proved. $\blacksquare$

\section{$A_{\infty}$-structure for the derived category of coherent sheaves}

\subsection{Rigid analytic space}

It will be helpful (although not necessary) for the reader
of this section
to be familiar with basic facts of non-archimedean analysis
(see [BGR]). For any smooth manifold $Y$ with integral affine structure
we will construct a sheaf ${\cal O}_Y$ of 
 ${\C}_{\varepsilon}$-algebras
 on $Y$.    Stalks ${\cal O}_{Y,y}$ of this sheaf are noetherian algebras, 
 and one can define the notion
 of coherent sheaves of   ${\cal O}_Y$-modules.
If $Y={\R}^n/{\Z}^n$ is the torus with the standard integral affine structure
then the category of coherent ${\cal O}_Y$-modules
 will be equivalent (by a non-archimedean version of GAGA) to the
category of coherent sheaves on an abelian variety over the
field ${\C}_{\varepsilon}$.

We start with the local picture. 
We denote by $v:{\C}_{\varepsilon}\to {\R}\cup \{+\infty\}$ a (non-discrete)
valuation defined by $v(\sum_{\lambda_1<\lambda_2<...}
c_ie^{-\lambda_i/\varepsilon})=\lambda_1$ if $c_1\ne 0$ and $v(0)=+\infty$.

\begin{dfn}
Let $U\subset {\R}^n$ be an open
subset of the standard vector space ${\R}^n$.
We define ${\cal O}_{{\R}^n}(U)$ as the vector space over ${\C}_{\varepsilon}$
consisting of formal Laurent series 
$$f=\sum_{k_1,...,k_n \in {\Z}^n}a_{k_1...k_n}z_1^{k_1}...z_n^{k_n},$$
where $z_1,...,z_n$ are formal variables, 
$a_{k_1...k_n}\in {\C}_{\varepsilon}$, and for any $(y_1,...,y_n)\in U$
we have: $lim_{\sum_i|k_i|\to \infty}(v(a_{k_1...k_n})+\sum_ik_iy_i)=+\infty$.

\end{dfn}

It follows from the definition that if $f\in {\cal O}_{{\R}^n}(U)$
and $(z_1,...,z_n)\in ({\C}_{\varepsilon}^{\ast})^n$ then the 
series $\sum_{k_1,...,k_n}a_{k_1...k_n}z_1^{k_1}...z_n^{k_n}$
converges in the adic topology as long as $(v(z_1),...,v(z_n))\in U$.

We introduce an action of the group $GL(n,{\Z})\ltimes {\R}^n$
on $({\R}^n,{\cal O}_{{\R}^n})$ such as follows:

a) $GL(n,{\Z})$ acts simultaneously by the linear change of coordinates and
linear transformation of indices $(k_1,...,k_n)$ in the series;

b) translations $(t_1,...,t_n)\in {\R}^n$ act on the coordinates
$(y_1,...,y_n)$ by the shift $(y_1,...,y_n)\mapsto (y_1+t_1,...,y_n+t_n)$, 
and on the series by the rescaling of coefficients
$$\sum_{k_1,...,k_n}a_{k_1...k_n}z_1^{k_1}...z_n^{k_n}\mapsto
\sum_{k_1,...,k_n}(a_{k_1...k_n}e^{-\sum_i  t_i k_i/\varepsilon})
z_1^{k_1}...z_n^{k_n}.$$

Using this action we define the sheaf ${\cal O}_{Y}$ for an arbitrary
smooth manifold $Y$ with integral affine structure.

We claim that there is a canonically associated
to $Y$ a rigid analytic space $Y^{an}$
defined over ${\C}_{\varepsilon}$. Here is the construction.
Let us consider a covering of $Y$ by open subsets $U_i$ such that
all non-empty intersections $U_{i_1i_2...i_k}:=U_{i_1}\cap...\cap U_{i_k}$
in some local affine coordinates are convex polyhedra whose faces have
rational slopes. Every  $U_{i_1i_2...i_k}$ can be identified with the
intersection of finitely many half-spaces, such that their
pre-images under the map $v^n:({\C}_{\varepsilon}^{\ast})^n\to {\R}^n$
are sets of the type $\{(z_1,\dots,z_n)|\,v(z_1^{k_1}...z_n^{k_n})\ge C\}$ 
for some
 rational $C>0$. It is known after Tate that such a system of inequalities
defines an affinoid domain (i.e. a local model for a rigid analytic space over 
${\C}_{\varepsilon})$.

\begin{dfn} We define $Y^{an}$ as
the rigid analytic space over ${\C}_{\varepsilon}$ obtained by
gluing the local data $(U_i,{\cal O}_{U_i})$ by means of the action of 
$GL(n,{\Z})\ltimes {\R}^n$.
\end{dfn}

It is easy to see that $Y^{an}$ is canonically defined, and that
the category $Coh(Y^{an})$ of coherent analytic sheaves on $Y^{an}$
(in the sense on analytic geometry) is equivalent to the category
of coherent ${\cal O}_Y$-modules (i.e. locally finitely generated
${\cal O}_Y$-modules).

To every algebraic variety ${\cal Y}$ over ${\C}_{\varepsilon}$
one can associate canonically a rigid analytic space
${\cal Y}^{an}$. If ${\cal Y}$ is projective then the category
$Coh({\cal Y}^{an})$ is equivalent to the category $Coh({\cal Y})$
of algebraic coherent sheaves on ${\cal Y}$ (GAGA theorem).

Assume that $Y={\R}^n/\Lambda$ is an
$n$-dimensional torus equipped with the standard
integral affine structure induced by ${\Z}^n\subset {\R}^n$,
and $\Lambda$ is a lattice commensurable with ${\Z}^n$.
The following result can be derived from [Mum].

\begin{prp} In the previous notation one has $Y^{an}\simeq {\cal Y}^{an}$
where ${\cal Y}$ is an abelian variety over ${\C}_{\varepsilon}$.

\end{prp}

Let us return to the picture of metric collapse 
in the case of abelian varieties.
Since the collapse was defined by rescaling of the lattice (see Section 2)
one can prove that ${\cal Y}$ is isomorphic to the original abelian
variety over ${\C}_{\varepsilon}$.
Therefore in the case of abelian varieties we have two equivalent
descriptions of the collapse: the one in terms of Riemannian geometry
and the one in terms of analytic non-archimedean geometry.

\begin{rmk} For the case of collapse with singular fibers,
 the rigid analytic space $Y^{an}$
constructed as above, seems to be a ``wrong'' one.
 First of all, it is not compact because $Y$ is not compact.  
 But there is also a more fundamental problem.
  It seems that $Y^{an}$ can not be embedded into a compact analytic space
   associated with a projective algebraic variety. 
   There are several indications that there exists another sheaf of algebras
    ${\cal O}'_Y$ which is (locally on $Y$) isomorphic to ${\cal O}_Y$,
     and the rigid analytic space $(Y^{an})'$ associated with $(Y,{\cal O}'_Y)$
     admits an algebraic compactification.
     In general, sheaves ${\cal O}'_Y$ which are twisted versions of 
     ${\cal O}_Y$
      are classified by the first non-abelian cohomology
       $H^1(Y,{\underline {Aut}}({\cal O}_Y))$ where 
       ${\underline {Aut}}({\cal O}_Y)$ is the sheaf of
       groups of automorphisms of ${\cal O}_Y$.   
       Thus, in the mirror symmetry for Calabi-Yau manifolds
        which are not abelian varieties, we expect
        a new ingredient, the cohomology class $[ {\cal O}'_Y]$.

        \end{rmk}

\subsection{$\A$-structure on coherent sheaves}

There is a sheaf of abelian groups ${Af}_Y$ on $Y$ given by locally
affine functions with integral slopes (such functions 
locally are given by $l=c+\sum_{1\le i\le n}m_iy_i$ where $m_i\in {\Z},
c\in {\R}$). There is a morphism of sheaves 
$exp:{Af}_Y\to {\cal O}_Y^{\ast}$ given by
$l\mapsto exp(l):=e^{-c/\varepsilon}\prod_{1\le i\le n}z_i^{m_i}$.

Let $(Y,g)$ be an AK-manifold (see Section 3.2).
We are going to define a characteristic class $[g]$ of the metric,
which will be an analog of the cohomology class of a
 K\"ahler form in complex geometry.
Let $Af_Y\otimes {\R}$ be a sheaf of all real-valued
locally affine functions. 
 For a cover by convex sets
$(U_i)_{i\in I}$ one can choose smooth functions $K_i$ such that
$g_{|U_i}=\partial^2 K_i$. Then $K_i-K_j\in Af_Y\otimes {\R}(U_i\cap U_j)$,
defines a 1-cocycle whose cohomology class we denote by
$[g]\in H^1(Y, Af_Y\otimes {\R})$.
If the {\it dual} affine structure (see Section 3)
is integral, we get a class 
$[g]$ in the subgroup $H^1(Y, Af_Y)/torsion\subset  H^1(Y, Af_Y\otimes {\R})$. 
We will call such classes {\it integral}.
In this case $exp([g])\in H^1(Y^{an}, {\cal O}_Y^{\ast})$
is the first Chern class of a line bundle on $Y^{an}$.
By analogy with the K\"ahler geometry we expect that this line
bundle is ample. In the case when $(Y,g)$ is a flat torus,
the ampleness can be proven directly (see [BL]).

From now on we assume that $[g]$ is integral.
Then by GAGA (see [Be], Prop. 3.14) the category of analytic coherent sheaves
on ${ Y}^{an}$ is equivalent to the category of algebraic
coherent sheaves on the corresponding algebraic projective variety
${\cal Y}$.

The sheaf ${\cal O}_Y$ admits a resolution $\widehat{\Omega}_Y^{\ast}$
by a soft sheaf of dg-algebras.
Locally, for a small open $U\subset Y$, sections of 
$\widehat{\Omega}_Y^{\ast}$ are given by sums
$\alpha=\sum_{i_1,...,i_n}c_{i_1...i_n}z_1^{i_1}...z_n^{i_n}$
where $c_{i_1...i_n}=\sum_jc_{j,i_1...i_n}
e^{-\lambda_{j,i_1...i_n}/\varepsilon},
c_{j,i_1...i_n}\in \Omega^{\ast}(U)$ with the same convergence conditions
as for the sheaf ${\cal O}_Y$.
Differential is given by the de Rham differential acting on the coefficients
$c_{j,i_1...i_n}$.

We define a dg-category ${\cal C}(Y)$ such as follows.
Objects are finite complexes of locally free ${\cal O}_Y$-modules
of finite rank. For any two such complexes $E_1$ and $E_2$
we define the space of morphisms as
$$Hom_{{\cal C}(Y)}(E_1,E_2)=\Gamma(Y,Hom_{{\cal O}_Y}(E_1,E_2)
\widehat{\otimes}_{{\cal O}_Y}\widehat{\Omega}_Y^{\ast}),$$
where we use the completed tensor product in the r.h.s. 
Differential and grading on the spaces of morphisms are
induced by those on $E_1,E_2, \widehat{\Omega}_Y^{\ast}$.
 We will treat ${\cal C}(Y)$ as an $\A$-pre-category in which all
 sequences of objects are transversal and there are no higher compositions
 except $m_1$ and $m_2$.

For a given projective algebraic variety $V$ over a field,
one can define canonically an equivalence class of $\A$-categories
$D^b_{\A}(V)$. It is obtained by the following enhancement
of the bounded derived category of
coherent sheaves on $V$. Objects of this $\A$-category are the same
as of the derived category of coherent sheaves.
In order to define
the space of morphisms between two objects, one replaces
them by arbitrary chosen acyclic resolutions 
by locally free sheaves (e.g. the Godement resolutions)
and then takes the global
sections of the space of morphisms between resolutions in the 
category of complexes of sheaves. In this way one obtains a dg-category.
 In the case of projective varieties over complex numbers,
there is an alternative construction in terms of complexes
of holomorphic vector bundles and Dolbeault forms.
Different choices of resolutions lead to
 $\A$-equivalent categories. We will denote the ($\A$-equivalence) 
 class
 of these categories by $D^b_{\A}(V)$.

By definition, spaces of morphisms of ${\cal C}(Y)$
are resolutions of the corresponding spaces of sheaves
of ${\cal O}_Y$-modules. Then using GAGA theorem from [Be], one concludes
that the following result holds.

\begin{prp} The category ${\cal C}(Y)$ is $\A$-equivalent
 to $D^b_{\A}({\cal Y})$, where
${\cal Y}$ is the projective algebraic variety corresponding
to the analytic space ${ Y}^{an}$ assigned to $Y$.

\end{prp}

\section{Homological mirror conjecture}

In the previous section we constructed an $\A$-category
${\cal C}(Y)$ which is $\A$-equivalent to the derived
category of coherent sheaves on a Calabi-Yau manifold over
the field ${\C}_{\varepsilon}$. In this section we are going to construct
a chain of $\A$-pre-categories and $\A$-equivalences 
(cf. with Section 6.6)
$${\cal C}_{unram}(Y)\mono {\cal C}_{unram,0}(Y) 
\hookleftarrow {\cal C}_{unram,0}^{tr}(Y)\hookleftarrow 
{\cal C}_{unram,0}^{tr,\Pi}(Y)
\leftarrow FO(X^{\vee})$$
and a functor $F:{\cal C}_{unram}(Y)\to {\cal C}(Y)$ which establishes
 an equivalence between ${\cal C}_{unram}(Y)$ and a full subcategory
 of ${\cal C}(Y)$. Recall that the Fukaya-Oh category $FO(X^{\vee})$,
  as defined in this paper, is also equivalent to a full
  subcategory of the Fukaya category $F(X^{\vee})$.
  Thus, we establish an $\A$-equivalence between full subcategories of
   the Fukaya category $F(X^{\vee})$ and of ${\cal C}(Y)$.

The  approach we are going to use is completely
similar to the one we used in the case of Morse
theory in Section 6. 
All the categories in our chain of $\A$-equivalence functors
 from above will have the same
 class of objects, i.e. the same as the Fukaya-Oh category.

\subsection{Mirror symmetry functor on objects}

Here we will define dg-category ${\cal C}_{unram}(Y)$ and the fully faithful 
 embedding $F$ of this category to ${\cal C}(Y)$. 
 
In the Appendix we will explain the conventional picture for the mirror
symmetry functor in case of complex numbers. There we will
 use a kind of Fourier-Mukai transform along fibers of the torus
fibration. The kernel of this transform is an analog
of Poincar\'e bundle. If one starts
with a local system on a Lagrangian section of $p^{\vee}:X^{\vee}\to Y$
then the transform makes from it a smooth bundle on $X$
with the connection which is flat in the anti-holomorphic directions.
In other words, one gets a holomorphic bundle on $X$.

These considerations cannot be literally repeated in the non-archimedean case.
We are going to construct the functor $F$ in the following way.
Let $(L,\rho)$ be an object of the category $FO(X^{\vee})$ such that
$rank(\rho)=1$ and the projection $L\to Y$ is one-to-one map. 
The manifold $L$ is locally
given by the graph of $df\,( mod\,(T_Y^{\Z})^{\vee})$, 
where $f$ is a smooth function on $Y$. To such an object we
assign a sheaf of rank one ${\cal O}_Y$-modules $F(L,\rho)$. 
For sufficiently small open $U\subset Y$ and chosen $f\in C^{\infty}(U)$
 the sheaf $F(L,\rho)_{|U}$ is identified with
${{\cal O}_Y}_{|U}$. 
Change $f\mapsto f+l$, where $l\in Af_Y(U)$
leads to the change of the trivialization of $F(L,\rho)_{|U}$ as
 $1_U\mapsto exp(l)\,1_U$
(here $1_U\in {\cal O}_Y(U)$ is the identity function). 
If $rank(\rho)$ is greater than one, we decompose $\rho_{|U}$ for small 
$U\subset Y$  into the sum of rank one local systems and then apply
the construction. Analogously, if the covering $L\to Y$ has more than 
one leaf, we apply the
previous construction to  each leaf of the covering and then
take the direct sum.

We will call $F$ the mirror symmetry functor on objects.
The category ${\cal C}_{unram}(Y)$ is defined as the dg-category whose class
 of objects
 is $Ob(FO(X^{\vee}))$, and the spaces of morphisms are
 $$Hom_{{\cal C}_{unram}(Y)}((L_1,\rho_1),(L_2,\rho_2)):=
 Hom_{{\cal C}(Y)}(F(L_1,\rho_1),F(L_2,\rho_2))$$
 
\subsection{Spectrum of a morphism and the semigroup}

Let $E_i=F(L_i,\rho_i), i=1,2$ be locally free  
${\cal O}_Y$-modules (i.e. vector bundles) corresponding to
objects $(L_i,\rho_i)\in FO(X^{\vee}), i=1,2$.
For any $\alpha\in Hom_{{\cal C}(Y)}(E_1,E_2)$ and a point $y\in Y$ 
we will define the spectrum of $\alpha$ at $y$ as a certain
(at most countable) discrete set of real numbers with finite multiplicities.

Let us assume first that $\rho_i, i=1,2$ are trivial rank one local
systems on $L_i, i=1,2$, and $L_i,i=1,2$ are unramified coverings of $Y$.
For a sufficiently small open set $U$ containing $y$ we 
can write in local coordinates
$L_i=graph(df_i)\,(mod\, (T_Y^{\Z})^{\vee}), i=1,2$ for smooth functions
$f_i:Y\to {\R}, i=1,2$.
Restriction to a small
open set $U$ of a morphism 
$\alpha\in Hom_{{\cal C}(Y)}(E_1,E_2)(U)=\widehat{\Omega}_Y^{\ast}(U)$
can be identified with the infinite series
$\alpha=\sum_{i_1,...,i_n}c_{i_1...i_n}z_1^{i_1}...z_n^{i_n}$, where
$c_{i_1...i_n}=\sum_jc_{j,i_1...i_n}
e^{-\lambda_{j,i_1...i_n}/\varepsilon}$ and 
$c_{j,i_1...i_n}\in {\Omega}_Y^{\ast}(U)$.

We define the
{\it spectrum of $\alpha$ at $y\in U$} as the set of real numbers
(with multiplicities)
$$Sp_y(\alpha)=\{-\lambda_{j,i_1...i_n}+
\sum_{1\le k\le n}i_ky_k+f_1(y)-f_2(y)\}\,,
$$
where the germ of $c_{j,i_1...i_n}$ at $y$ is not equal to zero.
One can check that $Sp_y(\alpha)$ is well-defined (i.e. does not depend
on the local trivialization),
and has the only limiting point at $s=-\infty$.

In the general case of higher rank local systems and Lagrangian
manifolds which are unramified coverings of $Y$, we decompose
$E_i, i=1,2$ locally near $y\in Y$ into the direct sum of
trivial rank one ${\cal O}_Y$-modules. The spectrum
of a morphism at the point $y$ is then defined as the union of the spectra of 
morphisms between corresponding line bundles.

\begin{rmk} One can use instead of the spectrum
an ${\R}$-filtration \\
$Hom_{{\cal C}(Y)}(E_1,E_2)^{\le s}$
on the space of morphisms.
It comes from the filtration on the stalks of sheaves
of morphisms 
$\underline{Hom}_{{\cal O}_Y}(E_1,E_2)
\widehat{\otimes}\widehat{\Omega}_Y^{\ast}$ (completed tensor product)
 defined by the condition 
$\{-\lambda_{j,i_1...i_n}+\sum_{1\le k\le n}i_ky_k+f_1(y)-f_2(y)\}\le s$.
It is easy to see that $\alpha$ belongs to $Hom_{{\cal C}(Y)}(E_1,E_2)^{\le s}$
iff for all $y\in Y$ one has $Sp_y(\alpha)\subset (-\infty,s]$.
\end{rmk}

Let us consider a subspace $Hom_{{\cal C}(Y)}^{alg}(E_1,E_2)\subset
Hom_{{\cal C}(Y)}(E_1,E_2)$ of algebraic morphisms.
It consists of finite sums (both in $z_i$
and $e^{-\lambda_{j,i_1...i_n}/\varepsilon}$). It is
dense in the space of all morphisms (analytic functions can be approximated
 by Laurent polynomials).
Moreover, the space $Hom_{{\cal C}(Y)}(E_1,E_2)$ coincides
with the completion of $Hom_{{\cal C}(Y)}^{alg}(E_1,E_2)$ with
respect to the ${\R}$-filtration introduced above.

There is a $1$-parameter semigroup $\phi^t, t\ge 0$ acting
on  $Hom_{{\cal C}(Y)}^{alg}(E_1,E_2)$. 
In local coordinates $\phi^t$ acts on the coefficients $c_{j,i_1...i_n}$
by moving them along the gradient flow of $f_1-f_2$.
In order to define it globally we need to describe the space
$Hom_{{\cal C}(Y)}^{alg}(E_1,E_2)$ in geometric terms.
It will be done below.

 Given two Lagrangian submanifolds
$L_i\subset X^{\vee}, i=1,2$ as above, a point $y\in Y$, two points
$x_i\in L_i, i=1,2$ such that $p^{\vee}(x_i)=y$,
we define a set $P(L_1,L_2,y)$
of homotopy classes of  paths $\gamma \in (p^{\vee})^{-1}(y)$
starting at $x_1$ and ending at $x_2$. Each homotopy class contains a unique
geodesic in the flat metric on the torus.
We define the space $P(L_1,L_2)=\sqcup_{y\in Y}P(L_1,L_2,y)$.
It carries an obvious topology such that the natural projection
$\pi:P(L_1,L_2)\to Y$ is an unramified covering with countable
fibers.
Using the symplectic form $\omega$ on $X^{\vee}$ we define a
closed $1$-form $\mu$ on $P(L_1,L_2)$ by the formula 
$\mu=\int_{\gamma}\omega$.
Locally on $Y$ we have: $L_i=df_i (mod\,(T_Y^{\Z})^{\vee}),
 i=1,2$ where $f_i:Y\to {\R}$
are smooth functions.
Then locally on $P(L_1,L_2)$ we have: $\mu=d(f_1-f_2+l)$,
where $l$ is a local section of the pullback of the sheaf $Aff_Y$. Clearly
the function $l$ is defined up the adding of a real constant.
Thus obtain an ${\R}$-torsor on $P(L_1,L_2)$.
Using the embedding ${\R}\to {\C}_{\varepsilon}^{\ast}$,
$\lambda\mapsto exp(\lambda/\varepsilon)$ we get a 
${\C}_{\varepsilon}^{\ast}$-torsor, which defines a local system
 ${\C}_{\varepsilon}^{tw}$
of $1$-dimensional ${\C}_{\varepsilon}$-modules over $P(L_1,L_2)$.
Fibers of  ${\C}_{\varepsilon}^{tw}$ carry natural
filtrations. Indeed, in a neighborhood of a point 
$(x_1,x_2,\gamma,y)\in P(L_1,L_2)$ we can choose a smooth function
$f=f_1-f_2+l$ such that $\mu=df$. It defines a local trivialization of 
${\C}_{\varepsilon}^{tw}$. In this trivialization the filtration
is defined for $h\in {\C}_{\varepsilon}$
by the condition $v(h)(y)+f(y)\le s, s\in {\R}$, where
$v$ is the valuation. We define a subsheaf  ${\C}_{\varepsilon}^{tw, alg}$
of ${\C}_{\varepsilon}^{tw}$ by the requirement that in a local
trivialization it is a subsheaf of finite sums of exponents.

Notice that there are natural projections $pr_i:P(L_1,L_2)\to L_i, i=1,2$.
Having local systems $\rho_i$ on $L_i, i=1,2$ we define local systems
$\widehat{\rho}_i, i=1,2$ on $P(L_1,L_2)$ as pullbacks with respect to
$pr_i, i=1,2$. 

On $P(L_1,L_2)$ we define a sheaf $\underline{Hom}^{alg}(E_1,E_2)$
($E_i, i=1,2$ were defined previously)
such as follows:
$\underline{Hom}^{alg}(E_1,E_2)={\C}_{\varepsilon}^{tw, alg}\otimes
(\widehat{\rho}_1)^{\ast}\otimes \widehat{\rho}_2\otimes 
\underline{\Omega}_{P(L_1,L_2)}^{\ast}$, 
where
$\underline{\Omega}_{P(L_1,L_2)}^{\ast}$ is the sheaf of differential forms.
We endow stalks of $\underline{Hom}^{alg}(E_1,E_2)$ with ${\R}$-filtrations
induced by the filtration on ${\C}_{\varepsilon}^{tw}$ and trivial
filtrations on the other tensor factors.

Let $\pi_{!}$ denotes the functor of direct image
with compact support. Then 
$\pi_{!}(\underline{Hom}^{alg}(E_1,E_2))=
\pi_!({\C}_{\varepsilon}^{tw}\otimes
\widehat{\rho}_1^{\ast}\otimes \widehat{\rho}_2)\otimes 
\underline{\Omega}_{Y}^{\ast}$, 
where the last tensor factor is the sheaf of de Rham differential
forms on $Y$.

We can identify ${\Z}^n$ with $H_1(T^n, {\Z})$, and the latter
group naturally acts on homotopy classes of paths $\gamma$.
On the other hand, the group ring of ${\Z}^n$ over ${\C}_{\varepsilon}$
can be identified with the ring of Laurent polynomials
${\C}_{\varepsilon}[z_1^{\pm 1},...,z_n^{\pm 1}]$. 
Let ${\C}_{\varepsilon}^{alg}\subset {\C}_{\varepsilon}$
be the subring of finite sums of exponents.
It is easy to see
that the structure of ${\C}_{\varepsilon}^{alg}[{\Z}^n]$-module
on the sections of
$\underline{Hom}^{alg}(E_1,E_2)$ corresponds to the structure
of ${\C}_{\varepsilon}^{alg}[z_1^{\pm 1},...,z_n^{\pm 1}]$-module
on its image under $\pi_!$.
 Using this observation one can prove
that 
$$Hom_{{\cal C}(Y)}^{alg}(E_1,E_2)\simeq
\Gamma(Y,\pi_{!}(\underline{Hom}^{alg}(E_1,E_2)))=
\Gamma_c(P(L_1,L_2),\underline{Hom}^{alg}(E_1,E_2)),$$
where the isomorphism is induced by the natural
morphism of sheaves 
$$\pi_{!}(\underline{Hom}^{alg}(E_1,E_2))\to
\underline{Hom}_{{\cal C}(Y)}^{alg}(E_1,E_2)\,.$$
 Here $\Gamma_c$ refers
to the functor of sections with compact support.

Using the metric on $Y$ we assign to the $1$-form $\mu$ a vector 
field $\xi$ on $P(L_1,L_2)$. Locally $\xi$ is the generator 
of the gradient flow of
$f_1-f_2+l$. It is not difficult to show that there is no trajectory
of the flow which goes to infinity for a finite time. Therefore
the vector field $\xi$  generates a $1$-parameter semigroup $\psi^t$ acting on $P(L_1,L_2)$.
The following result is easy to prove.

\begin{prp} The $1$-parameter semigroup $\psi^t$ decreases
the filtration on stalks of points
which do not belong to $L_1\cap L_2$. More precisely,
$$\psi^t(\underline{Hom}_p^{alg}(E_1,E_2)^s)\subset
\underline{Hom}_{\psi^t(p)}^{alg}(E_1,E_2)^{s-\int_0^t\mu},$$
where $p\in P(L_1,L_2)$ is an arbitrary point.

\end{prp}

Functor $\pi_!$ is compatible with the filtrations on the stalks
of sheaves $\underline{Hom}_{{\cal C}(Y)}^{alg}(E_1,E_2)$ and
$\underline{Hom}^{alg}(E_1,E_2)$.
It is easy to see that the completion of 
stalks of the former with respect to the filtration
induced from the one on $\underline{Hom}^{alg}(E_1,E_2)$ coincides with
$Hom_{{\cal C}(Y)}(E_1,E_2)$. Since the semigroup $\psi^t$ decreases
the filtration,
the semigroup $\phi^t$ extends continuously to the completion
with respect to the filtration.
Thus the following proposition holds.

\begin{prp} The action of $\phi^t$ extends continuously\\ 
from $Hom_{{\cal C}(Y)}^{alg}(E_1,E_2)$ to 
$Hom_{{\cal C}(Y)}(E_1,E_2)$.

\end{prp}

\subsection{Homological mirror symmetry for abelian varieties}

The approach and proofs are completely similar
to those from Section 6 (one has to
change the scalars to ${\bf C}_{\varepsilon}$ in Section 6).
We will omit the details. 

In the previous subsection we defined the semigroup
$\psi^t$ acting on the sections with compact support
$\Gamma_c(P(L_1,L_2),\underline{Hom}^{alg}(E_1,E_2))$.
This action corresponds to the action of the semigroup
$\phi^t$ on the space of morphisms $Hom_{{\cal C}(Y)}(E_1,E_2)$.
Notice that the sheaf 
$\underline{Hom}^{alg}(E_1,E_2)={\C}_{\varepsilon}^{tw, alg}\otimes
(\widehat{\rho}_1)^{\ast}\otimes \widehat{\rho}_2\otimes 
\underline{\Omega}_{P(L_1,L_2)}^{\ast}$ is a subsheaf of
${\C}_{\varepsilon}^{tw}\otimes
(\widehat{\rho}_1^{\ast}\otimes \widehat{\rho}_2)\otimes 
\underline{D}_{P(L_1,L_2)}^{\prime}$,
where
$\underline{D}_{P(L_1,L_2)}^{\prime}$ is the sheaf of 
distribution-valued differential forms on $P(L_1,L_2)$.
Sections of the latter sheaf carry the natural topology:
series in $exp(-\lambda_i/\varepsilon)$ with distributional
coefficients converge,
if they converge in the adic sense when paired with
a test differential form.
Similarly to the case of Morse theory (Section 6.5) one proves
the following result.

\begin{prp} For any $\beta \in 
\Gamma_c(P(L_1,L_2),\underline{Hom}^{alg}(E_1,E_2))$ there exists
limit  
$$\psi^{\infty}(\beta)=lim_{t\to +\infty}\psi^t(\beta)\in
\Gamma(P(L_1,L_2), {\C}_{\varepsilon}^{tw}\otimes
(\widehat{\rho}_1^{\ast}\otimes \widehat{\rho}_2)\otimes 
\underline{D}_{P(L_1,L_2)}^{\prime}).$$

\end{prp}

The limit is not difficult to describe in terms of the gradient flow
generating $\psi^t$. Using the fact that $\psi^t$ moves the spectrum
of a morphism to $-\infty$, one can prove similarly to the Section 6.5
that the limit $\psi^{\infty}(\beta)$ belongs to a finite-dimensional
${\C}_{\varepsilon}$-vector space generated by the distributions
corresponding the unstable manifolds 
$U_x\subset P(L_1,L_2), x \in L_1\cap L_2$. Clearly,
the map $\beta\mapsto \psi^{\infty}(\beta)$ extends to the completion
with respect to the filtration. It descends
to the map $\alpha \mapsto \phi^{\infty}(\alpha)$,
where $\alpha\in Hom_{{\cal C}(Y)}(E_1,E_2)$. The image
of $\phi^{\infty}$ belongs to the space isomorphic
to $Hom_{FO(X^{\vee})}((L_1,\rho_1),(L_2,\rho_2))$.

Then we repeat the arguments from Section 6.6. Namely,
we define an $\A$-pre-category ${\cal C}_{unram,0}(Y)$ 
in the same way as we defined the category
$DR_0(Y)$ in Section 6.6.
Objects of ${\cal C}_{unram,0}(Y)$ are the same as of ${\cal C}_{unram}(Y)$.
The spaces of morphisms
of ${\cal C}_{unram,0}(Y)$ are dg-modules over the dg-algebra 
${\C}_{\varepsilon}\widehat{\otimes}\Omega_0^{\ast}$,
where $\Omega_0^{\ast}$ is the dg-algebra of germs of differential
forms in the auxiliary parameter $\delta$
at $\delta=0\in {\R}_{\ge 0}$ (cf. Section 6.6). Compositions of morphisms
in ${\cal C}_{unram,0}(Y)$ are linear with respect to the
dg-module structure. The transversality conditions in
${\cal C}_{unram,0}(Y)$ and ${\cal C}_{unram}(Y)$ 
by definition are the same as in $FO(X^{\vee})$.
Thus we obtain an 
 $\A$-pre-category ${\cal C}_{unram,0}^{tr}(Y)$ which is 
$\A$-equivalent to ${\cal C}_{unram}(Y)$.

Using homological perturbation theory in the same way as in Section 6.6
(projectors and homotopies are now defined by means
of the semigroup $\phi^t$), we construct
the analog of the category $DR_0^{tr,\Pi}(Y)$. 
We denote this $\A$-pre-category
by ${\cal C}_{unrum,0}^{tr,\Pi}(Y)$. 
The spaces of morphisms of this category are completed tensor products of 
$\Omega_0^{\ast}$
 with finite-dimensional
${\C}_{\varepsilon}$-vector spaces, spanned by the
``smoothenings'' of the unstable currents $[U_x]$. These smoothenings
are defined by means of the operators $R_{\delta}$ in
the same way as in Section 6.5.
By definition, the category ${\cal C}_{unrum,0}^{tr,\Pi}(Y)$
 has the same transversality conditions
as the category $FO(X^{\vee})$. The spaces of morphisms are naturally
 quasi-isomorphic to the corresponding spaces of morphisms in
$FO(X^{\vee})$.
Repeating the arguments from Section 6.6 we will see
that the $\A$-structure on ${\cal C}_{unram,0}^{tr,\Pi}(Y)$ is equivalent
to the one on $FO(X^{\vee})$.
Indeed, we have a natural map from the space 
$Hom_{FO(X^{\vee})}((L_1,\rho_1),(L_2,\rho_2))$ 
to the space

$Hom_{{\cal C}_{unram,0}^{tr,\Pi}(Y)}(F(L_1,\rho_1),F(L_2,\rho_2))$.
Let us recall that all the categories ${\cal C}_{unram}(Y)$,
${\cal C}_{unram,0}(Y)$, ${\cal C}_{unram,0}^{tr}(Y)$ and
${\cal C}_{unrum,0}^{tr,\Pi}(Y)$ have the same class of objects.
Therefore we have the mirror symmetry functor
on objects $F:FO(X^{\vee})
\to {\cal C}_{unrum,0}^{tr,\Pi}(Y)$ (see Section 8.1).
On the other hand, the considerations above give rise
to the linear maps of the spaces of morphisms.
Let $\nu_{X_1,X_2}$ be this linear map
for two objects $X_i= (L_i,\rho_i), i=1,2$.
The proof of the following proposition is completely
similar to the corresponding one in Section 6.6.
(Since we work now over the field ${\bf C}_{\varepsilon}$,
one uses the formulas for $m_k, k\ge 1$ from the last
paragraph of Section 6.2).

\begin{prp} Let $E_i=F(X_i), 0\le i\le k, k\ge 1$
 be locally free rank one 
${\cal O}_Y$-modules (vector bundles) corresponding to
objects $X_i=(L_i,\rho_i)\in FO(X^{\vee}), 0\le i \le k$.
Then the formulas for 
$$m_k^{FO(X^{\vee})}:\otimes_{0\le i\le k}Hom(E_i,E_{i+1})
\to Hom(E_0,E_k)[2-k]$$ 
coincide (after the extension of scalars from ${\C}_{\varepsilon}$ to
 ${\C}_{\varepsilon}\widehat{\otimes}\Omega^{\ast}_0$) 
with the formulas
for 
$$m_k^{{\cal C}_{unram,0}^{tr,\Pi}(Y)}:
\otimes_{0\le i\le k}Hom(X_i,X_{i+1})\to Hom(X_0,X_k)[2-k]$$
when the spaces of morphisms are identified via the maps
$\nu(X_i,X_j)$.

\end{prp}

Therefore $\A$-pre-categories ${\cal C}_{unram,0}^{tr,\Pi}(Y)$ and
$FO(X^{\vee})$ are equivalent.
On the other hand, it follows from Sections 6.4, 6.6
 that ${\cal C}_{unram}(Y)$ and ${\cal C}_{unram,0}^{tr,\Pi}(Y)$
  are also equivalent.
  Finally, applying the functor $F$, we obtain the following theorem.

\begin{thm} The full subcategory $F({\cal C}_{unram}(Y))$ of $C(Y)$
is $\A$-equivalent to 
 $FO(X^{\vee})$. 
\end{thm}

This is our version of Homological Mirror Conjecture for abelian
varieties.

\begin{rmk} If we endow the torus 
$Y={\R}^n/{\Z}^n$ with a flat metric and consider
 only flat Lagrangian subtori in $X^{\vee}$ then all 
 higher compositions in the $\A$-pre-category $FO(X^{\vee})$
  can be written
  in terms of explicit ``truncated theta series'' analogous to those
  considered in [Ko] and [P1] in the case of elliptic curves.
\end{rmk}

\section{Appendix: constructions in the case of complex numbers}

In the previous section we considered algebraic and analytic varieties
over the complete local non-archimedean field ${\C}_{\varepsilon}$.
In this section we
explain our approach in the case of complex numbers (i.e. we will assume
that $\varepsilon$ is a fixed positive number). We should warn the reader that
it is not yet clear how to obtain rigorous proof of the Homological
Mirror Conjecture in this case.
In particular, it is not known how to prove convergence of
the series defining compositions in the Fukaya category.
Nevertheless we will discuss the complex case because
 the geometry is more transparent.
For example, one can construct the mirror symmetry functor on objects
by means of the real verison of Fourier-Mukai transform (see Section 9.1).
From the point of view of main body
of present paper, the Appendix can be treated as a geometric
motivation.
For this reason we will not stress that $X$ is an abelian variety,
but will be using our conjectures about the collapse, and the
assumption that the base $Y$ of the torus fibration
is a smooth manifold with integral affine structure
and K\"ahler potential. We will be also using the notation from Section 3.

\subsection{Mirror symmetry functor on objects over $\C$}

In the case of complex numbers the
mirror symmetry functor assigns a holomorphic vector bundle $F(L,\rho)$ 
on $X=X_{\epp}$
to a pair $(L,\rho)$, where $L\subset X^{\vee}$ is a Lagrangian submanifold,
such that the projection $p^{\vee}_{|L}:L\to Y$ is an unramified covering,
and $\rho$ is a 
local system on $L$. If $L$ is a section of $p^{\vee}$,
and $rank(\rho)=1$,
then $E=F(L,\rho)$ is a line bundle. In general, $E$ can be locally represented
as a sum $E\simeq \oplus_{\alpha \in A}E_{\alpha}$ where 
$A$ is the set of leaves
(i.e. connected components) of the covering $L\to Y$, and $E_{\alpha}$
is a holomorphic vector bundle of the rank equal to the rank of $\rho$ on 
the leaf $\alpha$.

 The following explicit construction
of the mirror symmetry functor on objects is not new, see e.g.
[AP].
We start with the remark that there is a canonical $U(1)$-bundle
on $X\times_Y X^{\vee}$ (Poincar\'e line bundle).
It will be denoted by ${P}$. It admits a  canonical
connection, which will be described below .
Let us fix $y\in Y$. Then $p^{-1}(y)\simeq T_{Y,y}/\epp T_{Y,y}^{\bf Z}$ and
$(p^{\vee})^{-1}(y)\simeq T_{Y,y}^{\ast}/(T_{Y,y}^{\bf Z})^{\vee}$.
We identify torus
$(p^{\vee})^{-1}(y)$ with the moduli space of  $U(1)$-local systems on
the torus $p^{-1}(y)$
trivialized over a point $0\in p^{-1}(y)$.  We define  $U(1)$-bundle    
$P$ to be the tautological
 bundle on   $X\times_Y X^{\vee}$ corresponding to this description.

In order to describe the connection on $P$
 let us consider the
fiberwise universal coverings $r:T_Y\to T_Y/T_Y^{\bf Z}$
and $r^{\vee}:T_Y^{\ast}\to T_Y^{\ast}/(T_Y^{\bf Z})^{\vee}$.
Then the pullback $\bar{P}$ of $P$ to $T_Y\times_YT_Y^{\ast}$
is canonically trivialized. Thus we can work in coordinates.
Let $y=(y_1,...,y_n)$ be coordinates on $Y$,
$x=(x_1,...,x_n)$ and $x^{\vee}=(x_1^{\vee},...,x_n^{\vee})$ 
be coordinates on the fibers of $T_Y\to Y$ and $T_Y^{\ast}\to Y$
respectively. Deck transformations   $x_j\mapsto x_j+\epp n_j,\,n_j\in \Z$ act on   $\bar{P}$
preserving the trivialization,
 and transformations    $x_j^{\vee}\mapsto  x_j^{\vee}+n_j^{\vee},\,n_j^{\vee}\in\Z$
 act  on   $\bar{P}$ by the multiplication    by
  $exp(2\pi i/\varepsilon \sum_j n_j^{\vee}x_j)$.

 Let $\nabla_0$ be the trivial connection
on  $\bar{P}$. We consider the connection $\bar{\nabla}$ on $\bar{P}$
which is given by the following formula

$$\bar{\nabla}=\nabla_0+2\pi i/\varepsilon \sum_{1\le j\le n}x_j^{\vee}dx_j.$$

\begin{lmm} The connection  $\bar{\nabla}$ gives rise to a 
connection on $P$.

\end{lmm}

{\it Proof.} Obviously,  connection
 $\bar{\nabla}$ does not change under the transformation
 $x_j\mapsto x_j+\epp n_j, n_j\in {\bf Z}$. The transformation
 $x_j^{\vee}\mapsto x_j^{\vee}+n_j^{\vee}, n_j^{\vee}\in {\bf Z}$ together with
 the gauge transformation of $\bar{\nabla}$
 by $h=exp(2\pi i/\varepsilon \sum_j n_j^{\vee}x_j)$ also preserves 
 $\bar{\nabla}$. This proves the Lemma.
$\blacksquare$

Let $(L,\rho)$ be as above. The mirror symmetry functor assigns
to it a holomorphic vector bundle $E=F(L,\rho)$ such that
(in coordinates) its fiber over a point $(y,x)$ is given by the
formula 
$E(y,x)=
\oplus_{\{x^{\vee}\in L,p^{\vee}(x^{\vee})=y\}}\rho(x^{\vee})
\otimes P(x,x^{\vee})$.
This vector bundle carries the induced connection
$\nabla_E$. In the case of unitary $\rho$ the bundle $E$ carries
also a natural hermitean metric.

\begin{prp} The  $(0,2)$-part of the curvature $curv(\nabla_E)$ is trivial.
In particular, $\nabla_E$ is a holomorphic connection.

\end{prp}

{\it Proof.} It follows from the fact that $L$ is Lagrangian.
Indeed, let us lift $L$ to $T_Y^{\ast}$.
Then locally in a neighborhood of a connected component of
$L$, one can find a smooth real function $f=f(y)$ such that
$L=df$. We can write the local equation for $L$:
$x_j^{\vee}=\partial f/\partial y_j,
1\le j\le n$. The connection $\nabla_E$ can be locally written as
 $\nabla_{E,0}+id_E\otimes (2\pi i/\epp \sum_j\partial f/\partial y_jdx_j)$,
 where $\nabla_{E,0}$ is the trivial flat connection on the
 vector bundle $E$.
Since the holomorphic coordinates on $T_Y$ are given by
$z_j=y_j+ix_j, i=\sqrt{-1}$, one sees that the $(0,2)$-part of the
curvature is equal to
$curv(\nabla_E)^{(0,2)}=
const\times
(\sum_{j,k}\partial^2f/\partial y_j\partial y_kd\bar{z_j}d\bar{z_k})=0$.
The Proposition is proved.
$\blacksquare$

\begin{dfn}
For any two holomorphic 
 vector bundles $E_1$ and $E_2$ on $X$, we define
$Hom_{Dolb}(E_1,E_2)=
\Omega^{0,\ast}(X,{Hom}(E_1,E_2))$.

\end{dfn}

  We consider the space
of Dolbeault differential forms with values in the vector bundle
$Hom(E_1,E_2)$ as a dg-algebra with respect to the 
$\bar{\partial}$-differential.
In this way one gets a structure of $A_{\infty}$-category
(in fact a dg-category) on the derived category of coherent sheaves
on $X$. One can show that this $\A$-structure is equivalent
to the one mentioned in the main text.

\subsection{Sectors in the space of Dolbeault forms}

Let $E_i=F(L_i,\rho_i), i=1,2$ be holomorphic vector bundles as above.
There is an analog of the dg-category ${\cal C}(Y)$ in the case
of complex numbers. We will denote it by ${\cal A}(Y)$. Objects
of ${\cal A}(Y)$ are holomorphic vector bundles on $X$ of the type
$E=F(L,\rho)$. Morphisms are sections of soft sheaves on $Y$.
Namely, we define the sheaf $\underline{Hom}_{{\cal A}(Y)}(E_1,E_2)$
on $Y$ as the direct image $p_{\ast}(\underline{Hom}_{Dolb}(E_1,E_2))$ 
(in the self-explained notation). Then $Hom_{{\cal A}(Y)}(E_1,E_2)$
are global sections of this sheaf.
This sheaf corresponds to the sheaf\\
$\underline{Hom}_{{\cal C}(Y)}(E_1,E_2)$
in the non-archimedean geometry.
Let us choose an open affine chart $U\subset Y, U\simeq {\R}^n$.
Then $\Gamma(U, \underline{Hom}_{{\cal A}(Y)}(E_1,E_2))$
contains a subsheaf of finite Fourier sums with respect
to the natural action of the torus\\
$T^n$ on
$\Gamma(U\times T^n, \underline{Hom}_{Dolb}(E_1,E_2))$.
Thus we have the sheaf $\underline{Hom}_{{\cal A}(Y)}^{alg}(E_1,E_2)$
which is an analog of the sheaf 
$\underline{Hom}_{{\cal C}(Y)}^{alg}(E_1,E_2)$ 
considered in the non-archimedean case.
Notice that there exists a natural homomorphism of sheaves
$j:p^{\ast}(\underline{\Omega}^{\ast}_Y)\to \underline{\Omega}^{0,\ast}_X$.
The image of $j$ consists of Dolbeault forms on $X$ which have
coefficients locally constant along fibers of $p$.
In local coordinates $j$ is given by the formula
$f_{i_1,...,i_n}(y)dy_{i_1}\wedge...\wedge dy_{i_n}\mapsto 
f_{i_1,...,i_n}(y)d\overline{z}_{i_1}\wedge...\wedge d \overline{z}_{i_n}$,
where $\overline{z}_k=y_k-\sqrt{-1}x_k, 1\le k\le n$.
It is easy to see that $j$ is compatible with the structure
of dg-algebras on de Rham and Dolbeault forms.
Thus for a pair of holomorphic vector bundles $E_1$ and $E_2$ on $X$
we have a canonical structure of dg-module over 
$\underline{\Omega}^{\ast}_Y$ on the sheaf
$\underline{Hom}_{{\cal A}(Y)}(E_1,E_2)$.
In the case when $E_i=F(L_i,\rho_i), i=1,2$ the subsheaf
$\underline{Hom}_{{\cal A}(Y)}^{alg}(E_1,E_2)$
is also a sheaf of dg-modules over $\underline{\Omega}^{\ast}_Y$.

As in the non-archimedean case there is a canonical decomposition
of the stalk $\underline{Hom}_{{\cal A}(Y)}^{alg}(E_1,E_2)_y, y\in Y$
into the direct sum of dg-modules of finite rank over
$\underline{\Omega}^{\ast}_Y$. Summands are labeled by
the homotopy classes $[\gamma]\in P(L_1,L_2,y)$ and
called {\it sectors}. We will denote them by
$\underline{Hom}_{{\cal A}(Y)}^{alg,[\gamma]}(E_1,E_2)_y$.
Informally, sectors correspond to ``Fourier components''
of Dolbeault forms in  $\underline{Hom}_{{\cal A}(Y)}^{alg}(E_1,E_2)_y$
in the direction of torus fibers.
Let us describe them more explicitly. For simplicity we will
assume that $\rho_i, i=1,2$ are rank one trivial local systems,
and $L_i, i=1,2$ intersect with each fiber of $p$ at exactly
one point. Then near $p^{-1}(y)$ we can write
$L_i=graph(df_i)\,( mod(T_Y^{\Z})^{\vee}), i=1,2$, where $f_i$
are germs at $y$ of smooth functions on $Y$.
From the description of the Poincar\'e bundle $P$ we deduce
that $\underline{Hom}_{{\cal A}(Y)}^{alg}(E_1,E_2)_y$
is canonically identified with the space of germs
of $\overline{\partial}$-forms near $T^n=p^{-1}(y)$,
endowed with the twisted differential
$\overline{\partial}^{\prime}\alpha=\overline{\partial}\alpha+
{i\over{\varepsilon}}
\sum\partial f/\partial y_id\overline{z}_i\wedge\alpha$,
where $f=f_1-f_2$.
Then the sector corresponding to a path $\gamma$ consists
of Dolbeault forms 
$\alpha=\sum_{i_1,...,i_n}exp(i\langle m,x\rangle/\varepsilon)
f_{i_1...i_n}(y)d\overline{z}_1\wedge...\wedge d\overline{z}_n$.
Here vector $m=m(\gamma)$ is the homotopy class of the loop
in $T^n=p^{-1}(y)$ which is the composition of three paths:

1) the path $[0,1]\to T^n, t\mapsto t(df_1)_y\, mod\,(T_Y^{\Z})^{\vee}$;

2) the path $\gamma$;

3) the path $[0,1]\to T^n, t\mapsto (1-t)(df_2)_y\, mod\,(T_Y^{\Z})^{\vee}$.

A choice of sector corresponds to the choice
of monomial $z_1^{i_1}...z_n^{i_n}$ in the non-archimedean case.
Homotopy classes of paths in non-archimedean approach correspond to
summands of Fourier series.
Locally each sector can be identified with the de Rham complex
on $Y$. Namely, to a form
$\alpha=\sum_{i_1,...,i_n}f_{i_1...i_n}(y)exp(\langle m, x \rangle)
d\overline{z}_{i_1}\wedge...\wedge d\overline{z}_{i_n}$
we assign the form $\alpha_m=
\sum_{i_1,...,i_n}f_{i_1...i_n}(y)exp({1\over{\varepsilon}}
\langle m, x \rangle)
dy_{i_1}\wedge...dy_{i_n}$, where $m=m(\gamma)$ defines the sector.
It is easy to see that
the differential $\overline{\partial}^{\prime}$ on
Dolbeault forms on $X$ corresponds to the de Rham differential
$d$ on $\Omega^{\ast}(Y)$.
In this way we obtain an isomorphism of complexes
$\underline{Hom}_{{\cal A}(Y)}^{alg, [\gamma]}(E_1,E_2)_y\simeq
\underline{\Omega}^{\ast}_{Y,y}\otimes {\C}$.

\begin{rmk} When $\varepsilon$ is not a fixed number,
but a parameter $\varepsilon \to 0$, 
 the coefficients $f_{i_1...i_n}(y)$ are asymptotic series
in $\varepsilon$ of the type
$f_{i_1...i_n}(y,\varepsilon)=\sum_{j\ge 1}
exp(-\lambda_j/\varepsilon)
f_{j,i_1...i_n}$ where $\lambda_j\in {\R}$ ,
$\lambda_1<...<\lambda_j<...$, and $\lambda_j\to +\infty$.

\end{rmk}

The set of exponents appearing in the expansion of $\alpha_m$ at $y$
corresponds to the spectrum $Sp_y(\alpha)$ considered
in the non-archimedean case.

\subsection{Semigroup $\varphi^t$}

Now we can define a semigroup
$\varphi^t:Hom_{Dolb}(E_1,E_2)\to Hom_{Dolb}(E_1,E_2),0\le t<+\infty$.
This is an analog of the semigroup $\phi^t$ in the non-archimedean case.
First, we identify the sector
$\underline{Hom}_{{\cal A}(Y)}^{alg, [\gamma]}(E_1,E_2)_y$
with $\underline{\Omega}^{\ast}(Y,Hom(\rho_1,\rho_2))_y$
as above.
Let us recall from the non-archimedean part, that to
the homotopy class of a path $\gamma$ we canonically associated
a closed $1$-form $\mu_{\gamma}=\int_{\gamma}\omega$, where $\omega$
is the symplectic form on $X^{\vee}$.
Using the Riemannian metric $g_Y$ on $Y$ we assign to
$\mu$ a vector field $\xi_{\gamma}$ on $Y$.
In a local trivialization it is given by
$grad((f_1-f_2+\langle m(\gamma),\cdot \rangle)/\varepsilon)$.
Then the infinitesimal action of $\varphi^t$ is defined
as the Lie derivative $Lie_{\xi_{\gamma}}$.
 Different Fourier components (sectors)
move on $Y$ with
with different speeds in different directions. Hence the picture is more 
complicated than in the case of Morse theory.

One can show that the generator $\Delta={d\over{dt}}|_{{t=0}}\varphi^t$
is a second order differential operator on $Hom_{Dolb}(E_1,E_2)$.
When $g_Y$ is a flat metric and $f_1=f_2=0$ one can find
the following explicit formula for $\Delta$:

$$\Delta=i\sum_j{\partial^2\over{\partial x_j\partial y_j}}-
{1\over{\varepsilon}}\sum_j{\partial^2\over{\partial x_j^2}}.$$ 

It seems plausible that there is an extension of
the semigroup $\varphi^t=e^{t\Delta}$, $t\ge 0$,
from $Hom_{{\cal A}(Y)}^{alg}(E_1,E_2)$ to
the whole space of morphisms $Hom_{{\cal A}(Y)}(E_1,E_2)$.
 Notice that $\Delta$ is not self-adjoint, and its real part is not elliptic.
 Nevertheless, we expect that the semigroup operator
 $\varphi^t$ converges as $t\to+\infty$ to a ``projector'' as 
 in the case of Morse theory.

\vspace{20mm}

{\bf References}

\vspace{5mm}

[AM] P.S.  Aspinwall, D. Morrison, String Theory  on K3 surfaces,
 hep-th/9404151.

\vspace{2mm}

[AP] D. Arinkin, A. Polishchuk, Fukaya category and Fourier transform,
math.AG/9811023.

 \vspace{2mm}

[Be] V. Berkovich, Spectral theory and analytic geometry over\\
non-Archimedean fields. Mathematical Surveys and Monographs, 33.
Amer. Math. Soc., 1990.

\vspace{2mm}

[BGR] S. Bosch, U. G\"unter, R. Remmert, Non-archimedean analysis.
Springer-Verlag, 1984.

\vspace{2mm}

[BL] S. Bosch, W. L\"utkebohmert, Degenerating abelian varieties,
Topology 30 (1991), no. 4, 653-698.

\vspace{2mm}

[BZ] J-M. Bismut, W. Zhang, An extension of a theorem by Cheeger and
M\"uller, Asterisque 205 (1992).

\vspace{2mm}

[CC] J. Cheeger, T.H. Colding, On the structure of spaces 
with Ricci curvature bounded below. I. J.
Differential Geom. 46 (1997), no. 3, 406-480.

\vspace{2mm}

[CG] J. Cheeger, M. Gromov, Collapsing Riemannian manifolds
while keeping their curvature bounded, Parts I and II,
 J. Diff. Geom., 23 (1986), 309-346
and 32 (1990), 269-298.

\vspace{2mm} 

[CY] S.-Y. Cheng, S.-T. Yau,
The real Monge-Ampere equation and affine flat structures,
Proceedings of the 1980 Beijing Symposium on Differential Geometry
and Differential Equations, vol.1, Science Press, Beijing 1982,
339-370.

\vspace{2mm} 

[De] P. Deligne, Local behavior of Hodge structures at infinity.
Mirror symmetry, II,
683-699, AMS/IP Stud. Adv. Math., 1, Amer. Math. Soc., Providence, RI, 1997

\vspace{2mm}

[Fu 1] K. Fukaya, Morse homotopy and its quantization, AMS/IP Studies 
in Adv. Math., 2:1 (1997), 409-440.

\vspace{2mm}

[Fu2] K. Fukaya, $\A$-category and Floer homologies. Proc. of GARC
Workshop on Geometry and Topology'93 (Seoul 1993), p. 1-102.

\vspace{2mm}

[FuO] K. Fukaya, Y.G. Oh, Zero-loop open string in the cotangent bundle
and Morse homotopy, Asian Journ. of Math., vol. 1(1998), p. 96-180).

\vspace{2mm}

[FuOOO] K. Fukaya, Y.G. Oh, H. Ohta, K. Ono, Lagrangian intersection
Floer theory-anomaly and obstruction. Preprint, 2000.

\vspace{2mm}

[G] M. Gromov, Metric structures for Riemannian manifolds,
 (J. Lafontaine and P. Pansu, editors), Birkh\"auser, 1999.

\vspace{2mm}

[Gaw] K. Gawedzki, Lectures on Conformal Field Theory, in Mathematical
Aspects of String Theory, AMS, 2000.

\vspace{2mm}

[Go] A. Goncharov,
Multiple zeta-values, Galois groups, and geometry of modular varieties,
math.AG/0005069 

\vspace{2mm}

[GS] V. Gugenheim, J. Stasheff, On perturbations and $\A$-structures,
Bul. Soc. Math. Belg. A38( 1987), 237-246.

\vspace{2mm}

[GW] M. Gross, P. Wilson, Large complex structure limits of K3 surfaces,
math.DG/0008018.

\vspace{2mm}

[H] N. Hitchin, The moduli space of special Lagrangian submanifolds,
math.DG/9711002.

\vspace{2mm}

[HL] F. Harvey, B. Lawson, Finite volume flows and Morse theory.
Preprint IHES 00/04, 2000.

\vspace{2mm}

[Ko] M. Kontsevich, Homological algebra of mirror symmetry. Proc. ICM Zuerich,
1994, alg-geom/9411018.

\vspace{2mm}

[KoS] M. Kontsevich, Y. Soibelman, Deformation theory, (book
in preparation).

\vspace{2mm}

[KoS1] M. Kontsevich, Y. Soibelman, Deformations of algebras over operads and
Deligne conjecture, math.QA/0001151.

\vspace{2mm}

[Le] N. Leung, Mirror symmetry without corrections, math.DG/0009235.

\vspace{2mm}

[LTY] B. Lian, A. Todorov, S-T. Yau,
Maximal Unipotent Monodromy for Complete Intersection CY Manifolds,
math.AG/0008061.

\vspace {2mm}

[M] M. Markl, Homotopy Algebras are  Homotopy Algebras, math.AT/9907138. 

\vspace{2mm}

[Me] S. Merkulov, Strong homotopy algebras of a K\"ahler manifold,
math. AG/9809172.

\vspace{2mm}

[Mo] D. Morrison, Compactifications of moduli spaces inspired by mirror
 symmetry, Journ\'ees de G\'eom\'etrie Alg\'ebrique d'Orsay (Juillet 1992), 
 Ast\'erisque, vol. 218, 1993, pp. 243-271, also alg-geom/9304007.

\vspace{2mm}

[Mum] D. Mumford, An analytic construction of degenerating abelian
varieties over complete rings. Compositio Math. 24:3 (1972), 239-272.

\vspace{2mm}

[P1] A. Polishchuk, $\A$-structures on an elliptic curve, 
math.AG/0001048.

\vspace{2mm}

[Se] P. Seidel, Vanishing cycles and mutations. math.SG/0007115.

\vspace{2mm}

[SYZ] A. Strominger, S-T. Yau, E. Zaslow, Mirror symmetry is $T$-duality,
hep-th/9606040.

\vspace{2mm}

[W1] E. Witten, Supersymmetry and Morse theory, J. Diff. Geom., v. 17 (1982),
661-692.

\vspace{5mm}

Addresses:

M.K.: IHES, 35 route de Chartres, F-91440, France

{maxim@ihes.fr}\\

Y.S.: Department of Mathematics, KSU, Manhattan, KS 66506, USA

{soibel@math.ksu.edu}

\end{document}